\documentclass[a4paper]{article}

\usepackage[english]{babel}
\usepackage[utf8x]{inputenc}
\usepackage[T1]{fontenc}

\usepackage[a4paper,top=3cm,bottom=2cm,left=3cm,right=3cm,marginparwidth=1.75cm]{geometry}

\usepackage{graphicx} 
\usepackage{subcaption} 

\usepackage{tabularray}
\usepackage{adjustbox}

\usepackage[hyphens]{url}
\usepackage{hyperref}
\usepackage{breakurl}

\usepackage{algorithm}
\usepackage{algpseudocode}

\algnewcommand{\Multcomment}[1]{\\\hfill \begin{minipage}[t]{0.8\textwidth}/* {\itshape #1} */\end{minipage}}

\algnewcommand{\algorithmicgoto}{\textbf{go to}}%
\algnewcommand{\Goto}[1]{\algorithmicgoto~\ref{#1}}%

\usepackage{enumitem} 
\usepackage{amsmath}
\usepackage{amssymb}
\usepackage{mathabx}
\usepackage[colorinlistoftodos]{todonotes}

\usepackage{color}
\usepackage{multirow}
\usepackage[round, sort,comma,authoryear]{natbib}
\usepackage{siunitx} 
\usepackage{verbatim}
\usepackage{tikz}

\usepackage{pifont}
\newcommand{\ev}{\protect\ding{51}}      
\newcommand{\best}{$\bigstar$}       
\newcommand{\rep}{$\blacklozenge$}        
\newcommand{\lead}{$\bigstar$}            

\title{Empirical Comparison of Unified Benders Cuts for Multi-Commodity Fixed-Charge Network Design}

\usepackage{authblk}   

\author[1]{Eric Larsen}
\author[2]{Jean-Fran\c{c}ois Cordeau}
\author[3]{Antonio Frangioni}
\author[1]{Emma Frejinger\thanks{Corresponding author: \texttt{emma.frejinger@umontreal.ca}}}

\affil[1]{Universit\'e de Montr\'eal, and CIRRELT, Montr\'eal, Canada}
\affil[2]{HEC Montr\'eal, and CIRRELT, Montr\'eal, Canada}
\affil[3]{Universit\`a di Pisa, Pisa, Italy}

\date{\today}

\begin{document}

\maketitle

\begin{abstract}
Among the many types of acceleration techniques designed to improve the performance of Benders decomposition, unified cut generation schemes have recently attracted a keen interest. Unified cuts aim for a better balance between the generation of optimality and feasibility cuts, while also providing a way to compare the strengths of different feasibility cuts. Our goal is to assess the experimental performance of a broad selection of unified and distinct Benders cuts in the context of the multi-commodity fixed-charge network design problem (MCFNDP). 

We express under a common mathematical structure and notation the construction of each unified or distinct Benders cut considered. We also explain how the generic formulations of the Benders cuts can be specialized to conform to the specifications of the MCFNDP. In addition, we suggest bespoke methods for comparing the performance of several solution methods when the benchmark is made up of heterogeneous problem instances. We report the results of a systematic empirical analysis comparing the performances of 50 Benders methods involving unified or distinct cuts in applications to a common testing bench made up of standardized MCFNDP instances. The analysis identifies a small number of leading Benders methods, namely those featuring the static Brandenberg–Stursberg cuts and the Hosseini–Turner l1-deepest cuts. In addition, we also report results obtained by using both Gurobi and CPLEX as the supporting solver to the SMS++ computation library.
\end{abstract}

\textbf{Keywords: Benders decomposition, unified cuts, network design problems, mixed integer linear programming}

\section{Introduction}
\label{sec:introduction}

Since its inception, \emph{Benders decomposition} \citep{Bend1962} has occupied a prominent position within the portfolio of decomposition methods aimed at reducing computation times and memory consumption in the context of large-scale mixed integer optimization. In particular, the application of Benders decomposition in association with the \emph{branch-and-cut} algorithm to deterministic mixed integer linear programs (MILPs) has been researched extensively and demonstrated to be highly effective. A prime application area for Benders decomposition is that of fixed-charge network design problems and their multiple variants \citep{Costa2005}. Successful applications also arise naturally in dual- or multi-stage stochastic programming with recourse when the probabilistic support is finite. Beyond these typical applications, there exists a broad range of extensions and refinements. They are discussed extensively in \cite{RahmEtAl2017}.

Among the numerous types of acceleration techniques introduced over the years to improve the performance of Benders decomposition, unified Benders cut generation has recently attracted a keen and sustained interest. Unified cuts aim to achieve a better balance between the generation of optimality and feasibility cuts while also offering a way to compare the strengths of different feasibility cuts. Our goal is to assess the experimental performance yielded by a broad selection of unified and distinct Benders cuts in the context of the multi-commodity fixed-charge network design problem (MCFNDP).

To promote a deeper understanding and a better usability of the techniques under consideration, we present them at first in their full generality within a canonical deterministic MILP setting. We select the following as our inceptive problem (IP):



\begin{align}
\min_{y, x}\; & f^T y + c^T x \label{IP:1} \\
\text{subject to }\; & G y + R x \leq q, \label{IP:2} \\
 & y \in \mathcal{Y}, \label{IP:3} \\
 & x \geq 0, \label{IP:4}
\end{align}
where $\mathcal{Y}$ is a bounded subset of $\mathbb{Z}^{n_y}$ or, without loss of generality, of ${\lbrace 0,1 \rbrace}^{n_y}$, for some $n_y$.

The contemporary \emph{branch-and-Benders-cut} algorithm tapping into the potential advantages offered by Benders decomposition combines the branch-and-cut algorithm of integer programming with a cutting plane algorithm in which so-called Benders cuts arise from the Benders decomposition. Before \cite{FiscEtAl2010}, \emph{feasibility} and \emph{optimality} cuts, respectively of the forms $\hat{\pi}^T(Gy -q) \leq 0$ and $\hat{\pi}^T(Gy - q) \leq \eta$, were usually generated separately. Starting with \cite{FiscEtAl2010} and continuing with \cite{Stur2019}, \cite{ConfWols2019}, \cite{HossTurn2021}, \cite{SeoEtAl2022}, and \cite{GlombEtAl2026}, several instances of the new breed known as \emph{unified} Benders cut have been introduced, with the aim of improving the computational performance of Benders decomposition in the solution of complex MILP problems.
Unified Benders cut for IP have the form $\hat{\pi}^T(Gy - q) - \hat{\pi}_0 \eta \leq 0$ where 
$(\hat{\pi}, \hat{\pi}_0)$ are simultaneously determined as the optimal solution of a bespoke \emph{cut-generating problem} (CGP), evaluated at the current candidate solution $\bar{y}$. These cuts are interpreted as either feasibility or optimality cuts depending on the value of $\hat{\pi}_0$ being zero or not.

\paragraph{Available evidence about the performance of unified cuts.}

To the best of our knowledge, the sum of comparative evidence currently available about the performances achieved by variants of the new unified cuts, variants of the preexisting \emph{distinct} cuts and direct application of branch-and-cut with MILP solvers is insufficient to issue a global judgment, even tentatively. The current evidence is fragmentary and disconnected: reported experiments involve (i) different subsets of the methods of interest, (ii) different sets of problem instances and (iii) different computational apparatuses. We are aware of the following sources: 
\begin{itemize}
\item \cite{FiscEtAl2010} report on the performance of Benders decomposition featuring their own proposal of a unified cut and that of a direct application of CPLEX when they are applied to small ad hoc instances of the multi-commodity, capacitated, fixed cost network design problem (MCFNDP) and to large standardized instances of the network expansion problem.
\item \cite{BranStur2021} report on the performance of their proposals of static and adaptive unified cuts when they are applied to large, actual instances of the network expansion problem.
\item \cite{SeoEtAl2022} report on their proposal of a unified cut, standard distinct cuts, Magnanti-Wong cuts and Fischetti et al. cuts. They are applied to ad hoc instances of MCFNDPs, to standardized instances of the communication network design problem, and to standardized instances of the network expansion problem.
\item \cite{HossTurn2021, HossTurn2025} report on their proposal of a unified cut adjoined with alternative underlying norms, standard distinct cuts and Fischetti et al. cuts. They are applied to standardized and to ad hoc instances of the capacitated facility location problem.
\item \cite{GlombEtAl2026} report on the performance of the following methods: (i) their proposals of a branch-and-Benders-cut algorithm equipped with its own unified Glomb et al. cut, (ii) alternative applications of Fischetti et al. cut, Brandenberg-Stursberg adaptive cut, Magnanti-Wong cut, (iii) a hybrid algorithm starting with the standard branch-and-Benders-cut adjoined with Fischetti et al. cut and finishing with the Glomb et al. algorithm, (iv) direct application of branch-and-cut with CPLEX. These methods are applied to the following sets of instances: ad hoc collection of instances from the MIPLIB repository; ad hoc creation of network design instances; network design instances from \cite{FiscEtAl2010}; ad hoc, randomly generated MILP problems. 
\end{itemize}

Table~\ref{tab:prior evidence} makes the structure of the prior evidence explicit: the \best~markers lie on the diagonal as each study reports its own, and typically its most elaborate cut as the winner, measured against a different and partial list of competitors on a different class of instances. Peeking ahead, our experiments revisit two conclusions that a reader might have drawn from the previous evidence: (i) Simplicity wins: the static Brandenberg–Stursberg cut and the Hosseini–Turner l1-deepest cuts lead, whereas the more sophisticated adaptive Brandenberg–Stursberg, Seo et al. and Glomb et al. cuts yield weaker performance on the MCFNDP. (ii) The direct branch-and-cut method is fastest overall; only Fischetti et al. (2010) and Glomb et al. (2026) included it as a competitor at all. The unified-cut literature may thus have underestimated the baseline it should outperform.

%

\begin{table}[htbp]
\centering
\small
\begin{adjustbox}{width=\textwidth}
\begin{tblr}{
  colspec = {l *{10}{c} l l},
  row{1} = {font=\bfseries},
}
\hline
Study & Dir. & Std & MW & Papa & Fisch & StatBS & AdaptBS & Seo & HT & Glomb & Instance class & Solver \\
\hline
Fischetti et al.\ (2010)        & \ev & \ev &        &        & \best  &        &        &        &        &        & MCFNDP (ad hoc); network expansion (std.)        & CPLEX \\
Brandenberg \& Stursberg (2021) &        &        &        &        &        & \ev & \best  &        &        &        & Network/power capacity expansion (real)          & Gurobi \\
Seo et al.\ (2022)              &        & \ev & \ev &        & \ev &        &        & \best  &        &        & MCFNDP (ad hoc); comm.\ network; net.\ expansion & CPLEX \\
Hosseini \& Turner (2021/25)    &        & \ev &        &        & \ev &        &        &        & \best  &        & Capacitated facility location (std.\ + ad hoc)   & CPLEX \\
Glomb et al.\ (2026)            & \ev &        & \ev &        & \ev &        & \ev &        &        & \best  & MIPLIB; net.\ design (incl.\ Fischetti et al.); rand.\ MILP & CPLEX \\
\hline
\textbf{This paper}             & \best  & \ev & \ev & \ev & \ev & \best  & \ev & \ev & \best  & \ev & MCFNDP Canad R (153, std., broad spectrum)        & Gurobi \& CPLEX \\
\hline
\end{tblr}
\end{adjustbox}
\caption{Available evidence on unified-cut performance.
\ev: method included in the study's comparison; \best: method the study reports as best. Each prior study evaluates its \emph{own} proposal against a different, partial set of competitors on a different instance class; none uses a common benchmark, a shared baseline of all families. On a common, broad MCFNDP benchmark this paper finds the \emph{static} Brandenberg--Stursberg cut and the Hosseini--Turner $l_1$-deepest cuts to be the leading Benders cuts, with direct branch-and-cut (Dir.) fastest overall; the adaptive Brandenberg--Stursberg, Seo et al. and Glomb et al. cuts -- each topping its own study -- do not lead here.}
\label{tab:prior evidence}
\end{table}

\paragraph{Experiments.}
The experiments we run compare the computational properties of 52 solution methods (50 Benders methods and 2 direct methods), each one combining an algorithmic family with a set of hyperparameters, applied to a testing bench consisting of instances of the linear, deterministic MCFNDP \citep{gendron1999multicommodity}. This stems from the recognition that network design problems occur in several critical areas such as supply chain management, transportation, and telecommunications and that instances of MCFNDP may serve as trial ground in stylized but near-realistic settings. Our test bench is made up of the Canad R set of standardized MCFNDP instances \citep{Canad}. It includes 153 valid generated problems featuring systematically differing numbers of nodes, arcs and commodities, as well as systematically differing magnitudes of fixed costs and arc capacities.

\paragraph{Contributions.}
We make the following contributions: (i) Expressing under the umbrella of a common mathematical structure and notation the construction of each formulation within an extensive selection of unified Benders cuts, as well as of some historically important non-unified Benders cuts. (ii) Offering a level of detail about the relevant CGPs, Benders cuts, and variants of the branch-and-Benders-cut algorithm that is sufficient to conduct applications. (iii) Detailing how the generic formulation of IP must be specialized to conform to the specifications of the MCFNDP. (iv) Supplying element-explicit descriptions in connection with applications to MCFNDP for every CGP and Benders cut considered. (v) Suggesting bespoke methods for comparing the performance of several computation methods in experiments where (va) existence and availability of some raw measurements may depend on experimental outcomes and/or (vb) magnitude of raw measurements may depend on the widely varying sizes and complexity of the problems included in the experiments. (vi) Performing a systematic empirical analysis that evaluates and compares the performances yielded by a broad selection of high-interest computational methods featuring unified and/or non-unified Benders cuts, in applications to a common testing bench made up of a broad spectrum of MCFNDP instances. The analysis identifies a small number of leading Benders methods among the 50 that are evaluated. (vii) Running parallel sets of experiments where Gurobi or CPLEX play the role of supporting solver so that their computational properties may be compared. (viii) Illustrating the effectiveness of the SMS++ library \citep{SMS++} performing the backend computations.

\paragraph{Structure and content.}
The remainder of the paper is structured as follows: Section~\ref{sec:Benders algo} overviews the Benders decomposition and the branch-and-Benders-cut algorithm for the purpose of an empirical investigation. It describes the construction of the standard feasibility and optimality cuts, the historically important and empirically successful Magnanti-Wong \citep{MagnWong1981} and Papadakos \citep{Papa2008} non-unified optimality cuts and specifies the relevant CGPs. Section~\ref{sec:unified cuts} details the unified Benders cuts introduced by \cite{FiscEtAl2010, BranStur2021, SeoEtAl2022, HossTurn2025, GlombEtAl2026}. For each one, it defines the particular CGP and Benders cut. Section~\ref{sec:MCFNDP} presents the MCFNDP underlying the comparative experiments. It includes an element-explicit statement of the MCFNDP and a matrix representation of the MCFNDP as a specialization of IP. Section~\ref{sec:description of experiments} details the set of experiments that are conducted. It describes the algorithmic families and the hyperparameters making up the computational methods, the set of standardized MCFNDP instances used as a testing bench, the computational apparatus. Section~\ref{sec:measuring performance} describes how the computational performances of all methods are assessed and compared with bespoke methods and with computation profiles. Section~\ref{sec:experimental results} details the results of the experiments. It identifies leading methods and discusses the computational characteristics of Gurobi and CPLEX. Section~\ref{sec:conclusion} summarizes our initial goals and our findings. It outlines envisioned extensions. The appendix supplies element-explicit expansions for every CGP and Benders cut included in the experiments.

\section{Benders decomposition with distinct feasibility and optimality cuts}
\label{sec:Benders algo}

Expressing IP with the nested formulation
\begin{equation}
\min_{y \in \mathcal{Y}} \lbrace f^T y + \min_{x \geq 0} \lbrace c^T x \mid Rx \leq q - Gy \rbrace \rbrace \label{eq:telesc formul}
\end{equation}
 hints of a possible separation of concerns between a delegating outer problem and an inner, slave problem. Indeed, starting from a restriction of the outer problem, a modern algorithm featuring Benders decomposition iterates on the following cycle within a branch-and-cut enumerative process that progressively explores $\mathcal{Y}$: (i) Calculate a tentative solution for the outer problem. (ii) Check the validity of the tentative solution within the inner problem. If it is invalid, define a half-space that is consistent with (\ref{IP:2}), (\ref{IP:3}) but excludes the current solution and insert it as an additional constraint in the restricted outer problem. Otherwise, update the incumbent solution (the value of $y$ producing the lowest higher bound) if appropriate. When the enumerative process has finished exploring $\mathcal{Y}$, the incumbent solution coincides with that of IP.
 
In (\ref{eq:telesc formul}) the inner problem and its dual are known as the Benders \emph{primal subproblem} (PSP) and the Benders \emph{dual subproblem} (DSP), respectively:
\paragraph{PSP:}
\begin{align}
\min_{x}\; & c^T x \label{PSP:1}\\
\text{subject to }\; & Rx \leq q - Gy, \label{PSP:2}\\
 & x \geq 0. \label{PSP:3}
\end{align}
\paragraph{DSP:}
\begin{align}
\max_{\hat{\pi}}\; & (Gy - q)^T \hat{\pi} \label{DSP:1}\\
\text{subject to }\; & -R^T\hat{\pi} \leq c, \label{DSP:2}\\
 & \hat{\pi} \geq 0. \label{DSP:3}
\end{align}
Defining PSP's value function as $\widehat{\textit{V}}(y) :\equiv \min_{x \geq 0} \lbrace c^T x \mid Rx \leq q - Gy \rbrace$, $\widehat{\textit{V}}$'s epigraph is
\begin{equation}
epi(\widehat{\textit{V}}) = \lbrace (y, \eta) \in conv(\mathcal{Y}) \times \mathbb{R} \mid \exists \,x \geq 0 :  Rx \leq q - Gy, c^T x \leq \eta \rbrace,
\end{equation}
and IP can be expressed as $\min \lbrace f^T y + \eta \mid (y, \eta) \in epi(\widehat{\textit{V}}) \cap (\mathcal{Y} \times \mathbb{R})  \rbrace.$ Now, by strong duality $\widehat{\textit{V}}(y) = \max_{\hat{\pi} \geq 0} \lbrace (Gy - q)^T \hat{\pi} \mid -R^T \hat{\pi} \leq c \rbrace,$ whence
\begin{equation}
epi(\widehat{\textit{V}}) = \lbrace(y, \eta) \in conv(\mathcal{Y}) \times \mathbb{R} \mid 
(Gy - q)^T \hat{\pi}^r \leq 0, \forall \hat{\pi}^r \in \widehat{\Pi}^R,  (Gy - q)^T \hat{\pi}^e \leq \eta, \forall \hat{\pi}^e \in \widehat{\Pi}^E \rbrace,
\end{equation}
where $\widehat{\Pi}^E$ is the set of extreme points of $\widehat{\Pi} :\equiv \lbrace \hat{\pi} \geq 0, -R^T \hat{\pi} \leq c \rbrace$ and $\widehat{\Pi}^R$ is the set of extreme rays of $\lbrace \hat{\pi} \geq 0, -R^T \hat{\pi} \leq 0 \rbrace$, i.e., the recession cone of $\widehat{\Pi}$.
Thus, IP can also be expressed as 
\begin{align}
\min_{y, \eta}\; & f^T y + \eta \label{MP:1}\\
\text{subject to }\; & (Gy - q)^T \hat{\pi}^e \leq \eta, \forall \hat{\pi}^e \in \widehat{\Pi}^E, \label{MP:2}\\
& (Gy - q)^T \hat{\pi}^r \leq 0, \forall \hat{\pi}^r \in \widehat{\Pi}^R, \label{MP:3}\\
& y \in \mathcal{Y}. \label{MP:4}
\end{align}
This is the so-called Benders \emph{master problem} (MP). Constraints (\ref{MP:2}) and (\ref{MP:3}) are known as the standard \emph{optimality} and \emph{feasibility} cuts, respectively.

The modern algorithm combining Benders decomposition with branch-and-cut, known as \emph{branch-and-Benders-cut} \citep{FortzPoss2009,Gendron2016}, starts by creating (i) an enumerative branch-and-cut process over $\mathcal{Y}$ equipped with a single search tree and (ii) a restricted master problem (RMP), where constraints sets (\ref{MP:2}) and (\ref{MP:3}) are typically removed entirely and by (iii) introducing a safe lower bound on $\eta$. The algorithm then progressively outer-approximates $epi(\widehat{\textit{V}})$ by launching the execution of the following sequence of steps through a callback mechanism each time a new candidate solution $(\bar{y}, \bar{\eta})$ to RMP where $\bar{y}$ is integral is found by the enumerative process:

\paragraph{Benders iteration.}
(i) Verify if PSP is feasible or not (i.e., DSP is bounded or not) at the current candidate value $\bar{y}$. (ia) If PSP is infeasible (i.e., DSP is unbounded), select a ray in the recession cone of $\widehat{\Pi}$, build the corresponding feasibility cut and add it as a lazy constraint to RMP. (ib) If PSP is feasible (i.e., DSP is bounded) but $\widehat{\textit{V}}(\bar{y}) > \bar{\eta}$, build the optimality cut from the solution of DSP and add it as a lazy constraint to RMP.
(ic) Otherwise, the current candidate value $\bar{y}$ is admissible for IP and
$f^T\bar{y} + \widehat{\textit{V}}(\bar{y})$ is a valid upper bound for the optimal value of
IP. Update the incumbent solution with the candidate value if appropriate. (ii) If the branch-and-cut enumerative process over $\mathcal{Y}$ is not completed yet, pursue it with the updated RMP until a new integral solution is found. (iia) If and when this occurs, return to (i). (iib) Otherwise, the incumbent solution is chosen as the final solution to IP.

\subsection{Magnanti-Wong optimality cut}
\label{sec:Magnanti-Wong optim cut}

\cite{MagnWong1981} define the following CGP:
\paragraph{CGP-M:}
\begin{align}
\max_{\hat{\pi}}\; &  (Gy^* - q)^T \hat{\pi} \label{CGP-M:1}\\
\text{subject to }\; & (G\bar{y} - q)^T \hat{\pi} = \widehat{\textit{V}}(\bar{y}), \label{CGP-M:2}\\
& -R^T \hat{\pi} \leq c, \label{CGP-M:3}\\
 & \hat{\pi} \geq 0, \label{CGP-M:4}
\end{align}
where the reference point $y^* \in relint\, conv(\mathcal{Y})$ and $\bar{y}$ is the current RMP solution. The resulting Magnanti-Wong optimality cut is given by $(Gy - q)^T \hat{\pi} \leq \eta$, where $\hat{\pi}$ is the current solution to CGP-M. If $y^* \in relint\, conv(\mathcal{Y})$, the cut is Pareto-optimal: Given $\hat{\pi}^0$ yielded by CGP-M, there is no other $\hat{\pi}^1$ such that $(Gy - q)^T \hat{\pi}^1 \geq (Gy - q)^T \hat{\pi}^0$ for all admissible $y$ and $(Gy - q)^T \hat{\pi}^1 > (Gy - q)^T \hat{\pi}^0$ for at least one admissible $y$.

\paragraph{Benders iteration.}
Similar to that of Section~\ref{sec:Benders algo}, except that the Magnanti-Wong optimality cut is added at step (ib), in addition to the standard optimality cut.

\subsection{Papadakos optimality cut}
\label{sec:Papadakos optim cut}

\cite{Papa2008} defines instead this CGP:
\paragraph{CGP-P:}
\begin{align}
\max_{\hat{\pi}}\; &  (Gy^* - q)^T \hat{\pi} \label{CGP-P:1}\\
\text{subject to }\; & -R^T \hat{\pi} \leq c, \label{CGP-P:2}\\
 & \hat{\pi} \geq 0, \label{CGP-P:3}
\end{align}
where $y^* \in relint\, conv(\mathcal{Y})$. The Papadakos optimality cut is given by $(Gy - q)^T \hat{\pi} \leq \eta$, where $\hat{\pi}$ is the current solution to CGP-P.
If $y^* \in relint\, conv(\mathcal{Y})$ the cut is Pareto-optimal. Here also, application of a standard feasibility cut is required when PSP is infeasible. While CGP-P is generally faster to solve than CGP-M, the Papadakos cuts are decoupled from the current solution of the RMP due to the absence of constraint (\ref{CGP-M:2}), unless the reference point $y^*$ is consistently updated as suggested in Papadakos (2008).

\paragraph{Benders iteration.}
Similar to that of Section~\ref{sec:Benders algo}, except that the Papadakos optimality cut is added at step (ib), in addition to the standard optimality cut.

\section{Unified Benders cuts}
\label{sec:unified cuts}

A unified Benders cut for IP has the form $\hat{\pi}^T(Gy - q) - \hat{\pi}_0 \eta \leq 0$ where 
$(\hat{\pi}, \hat{\pi}_0)$ are simultaneously determined as the optimal solution of a CGP evaluated at the current candidate value $\bar{y}$ and added as a lazy constraint to the RMP, provided a number of conditions regarding the feasibility, boundedness and optimal value of the CGP are satisfied. These cuts are interpreted as either feasibility or optimality cuts depending on the value taken by $\hat{\pi}_0$ being zero or not. As far as we know, original unified cuts appear in \cite{FiscEtAl2010}, \cite{Stur2019}, \cite{ConfWols2019}, \cite{HossTurn2021}, \cite{SeoEtAl2022} and \cite{GlombEtAl2026}.

\subsection{Unified cuts of Fischetti et al., Brandenberg-Stursberg, Conforti-Wolsey}
\label{sec:unified cuts of FBC}

\cite{FiscEtAl2010, BranStur2021, ConfWols2019} all build on the seminal work of \cite{CornLema2006} to define their respective versions of a unified Benders cut. In contrast with Section~\ref{sec:Benders algo}, they treat the two defining conditions of $epi(\widehat{\textit{V}})$ jointly rather than separately. By Farkas' lemma, $(\bar{y}, \bar{\eta}) \in epi(\widehat{\textit{V}})$ iff the following \emph{certificate-generating problem} has a bounded optimal value (equal to zero):
\begin{align}
\max_{\hat{\pi}, \hat{\pi}_0}\; & \hat{\pi}^T(G\bar{y} - q) - \hat{\pi}_0 \bar{\eta} \\
\text{subject to }\; & (\hat{\pi}, \hat{\pi}_0) \in \widehat{\Pi} :\equiv \lbrace \hat{\pi} \geq 0, \hat{\pi}_0 \geq 0 \mid \hat{\pi}^T R + \hat{\pi}_0 c^T \geq 0 \rbrace.
\end{align}
If $(\bar{y}, \bar{\eta}) \notin epi(\widehat{\textit{V}}), \exists \, (\hat{\pi}, \hat{\pi}_0) \in \widehat{\Pi}$ (a Farkas \emph{certificate}) $\mid \hat{\pi}^T(G\bar{y} - q) - \hat{\pi}_0 \bar{\eta} > 0$ whence $(\hat{y}, \hat{\eta})$ breaches the valid constraint $\hat{\pi}^T(Gy - q) - \hat{\pi}_0 \eta \leq 0$. When added to MP or RMP, the latter is known as a \emph{unified cut}.  Hence, $epi(\widehat{\textit{V}}) = \lbrace (y, \eta) : (y, \eta) \in \widehat{\mathcal{H}}(\hat{\pi}, \hat{\pi}_0), \forall \, (\hat{\pi}, \hat{\pi}_0) \in \widehat{\Pi} \rbrace$, where $\widehat{\mathcal{H}}(\hat{\pi}, \hat{\pi}_0) :\equiv \lbrace (y, \eta) : \hat{\pi}^T(Gy -q) - \hat{\pi}_0 \eta \leq 0 \rbrace$, and MP may be defined based on unified cuts as
\begin{align}
\min_{y, \eta}\; & f^T y + \eta \\
\text{subject to }\; & (y, \eta) \in \widehat{\mathcal{H}}(\hat{\pi}, \hat{\pi}_0), \forall (\hat{\pi}, \hat{\pi}_0) \in \widehat{\Pi}, \\
 & \eta \in \mathbb{R}, y \in \mathcal{Y}.
\end{align}

\cite{FiscEtAl2010, BranStur2021, ConfWols2019} differ with respect to the manner in which, at Benders iterations of the branch-and-Benders-cut algorithm, they select the Farkas certificates $(\hat{\pi}, \hat{\pi}_0) \in \widehat{\Pi}$ defining the unified cuts that are added to the RMP. Two CGPs may be distinguished here: CGP-FBC and CGP-FBC-invert. Given values assigned to $\tilde{\omega}, \tilde{\omega}_0$, CGP-FBC selects particular $(\hat{\pi}, \hat{\pi}_0) \in \widehat{\Pi}$. CGP-FBC-Invert reverses the roles of the objective and normalization constraint appearing in CGP-FBC \citep{CornLema2006}. Although CGP-FBC and CGP-FBC-Invert have nearly identical mathematical properties, their computational stability and speed may differ. In both cases, $(\bar{y}, \bar{\eta})$ is the current integral solution of RMP found by the ongoing enumerative branch-and-cut process over $\mathcal{Y}$.

\paragraph{CGP-FBC:}
\begin{align}
  \max_{\hat{\pi}, \hat{\pi}_0}\; & \tilde{\omega}^T \hat{\pi} + \tilde{\omega}_0 \hat{\pi}_0 \label{CGP-FBC:1}\\
  \text{subject to }\; & R^T \hat{\pi} + c \hat{\pi}_0 \geq 0, \label{CGP-FBC:2}\\
  & (q - G\bar{y})^T \hat{\pi} + \bar{\eta} \hat{\pi}_0 \leq -1, \label{CGP-FBC:3}\\
  &\hat{\pi}, \hat{\pi}_0 \geq 0. \label{CGP-FBC:4}
\end{align}
The feasible domain of CGP-FBC, when it is not empty, contains all normalized Farkas certificates from which a unified cut based on Farkas certificates may be constructed. When CGP-FBC is infeasible, no cut is required and the current solution $(\bar{y}, \bar{\eta})$ is admissible for IP. The feasible domain of CGP-FBC is called the \emph{relaxed alternative polyhedron} \citep{BranStur2021} or \emph{alternative polyhedron} \citep{GleeRyan1990, FiscEtAl2010} if the normalization (\ref{CGP-FBC:3}) appears instead as an equality. Both formulations lead to the same solution.

\paragraph{CGP-FBC-Invert:} \label{para:CGP-FBC-Invert}
\begin{align}
  \max_{\hat{\pi}, \hat{\pi}_0}\; & (G\bar{y} - q)^T \hat{\pi} - \bar{\eta} \hat{\pi}_0 \label{CGP-FBC-Invert:1}\\ 
  \text{subject to } & R^T \hat{\pi} + c\, \hat{\pi}_0 \geq 0, \label{CGP-FBC-Invert:2}\\
  &\tilde{\omega}^T \hat{\pi} + \tilde{\omega}_0 \hat{\pi}_0 = -1, \label{CGP-FBC-Invert:3}\\
  &\hat{\pi}, \hat{\pi}_0 \geq 0. \label{CGP-FBC-Invert:4}
\end{align}
When the optimal value of CGP-FBC-Invert is zero, no cut is required and the current solution $(\bar{y}, \bar{\eta})$ is admissible for IP.

\paragraph{Benders cuts.}
The unified cut resulting from either CGP-FBC or CGP-FBC-Invert is as follows:
\begin{equation}
(q - Gy)^T \hat{\pi} + \eta\hat{\pi}_0 \geq 0. \label{CGP-FBC_cut:1}
\end{equation}
Whenever $\hat{\pi}_0 = 0$, (\ref{CGP-FBC_cut:1}) turns out to be a feasibility cut. Otherwise, it is an optimality cut.
Particular values assigned to $\tilde{\omega}, \tilde{\omega}_0$ give rise to the cuts described in \cite{FiscEtAl2010, BranStur2021, ConfWols2019}.

\cite{FiscEtAl2010} selects a vertex of the \emph{alternative polyhedron} and therefore a particular \emph{minimal infeasibility system} (MIS) among the rows of the system  $Rx \leq q - G\bar{y}, c^T x \leq \eta$, by setting $\tilde{\omega}_0 = -1$ and
\begin{equation}
  (\tilde{\omega}_i) :=
  \begin{cases}
    0, & \text{if } G_{ij} = 0, \forall j\\
    -1, & \text{otherwise}.
  \end{cases} 
\end{equation}
This follows the intuition that stronger cuts and faster computation should result from concentrating cuts on sets of active sources of infeasibility. Nevertheless, the resulting cuts are not guaranteed to support $epi(\widehat{\textit{V}})$.

\cite{BranStur2021} demonstrate that setting $(\tilde{\omega}, \tilde{\omega}_0) = (G \omega, -\omega_0)$ for some $(\omega, \omega_0)$ selects a vertex of the \emph{reverse polar set} (see \cite{Bala1979}) constructed from $epi(\widehat{\textit{V}})$ (and, as a result, selects a vertex of the alternative polyhedron as well). The resulting cut is guaranteed to support $epi(\widehat{\textit{V}})$, and, if the solution $(\hat{\pi}, \hat{\pi}_0)$ is unique, it either supports a facet of $epi(\widehat{\textit{V}})$ or possibly contains $epi(\widehat{\textit{V}})$ if $epi(\widehat{\textit{V}})$ is not full-dimensional. Moreover, if $(\omega, \omega_0) \in relint(conv(epi(\widehat{\textit{V}}) \cap (\mathcal{Y} \times \mathbb{R})-(\bar{y}, \bar{\eta}))$ and $\hat{\pi}_0 > 0$, then the resulting optimality cut is Pareto-optimal. \cite{BranStur2021} distinguish two particular cases: (i) the \emph{static} cut arising from $(\tilde{\omega}, \tilde{\omega}_0) = (G \mathbf{1}, -1)$, (ii) the \emph{adaptive} cut arising from $(\tilde{\omega}, \tilde{\omega}_0) = (G (\tilde{y} - \bar{y}), -(\tilde{\eta} - \bar{\eta}))$, where $(\bar{y}, \bar{\eta})$ is the current solution of RMP and $(\tilde{y}, \tilde{\eta})$ is the current incumbent solution to IP calculated with RMP. The latter ensures that any optimality cut that is generated is Pareto-optimal. \cite{ConfWols2019} describe a similar cut.

\paragraph{Benders iteration.}
(i) Verify that CGP-FBC is feasible or that CGP-FBC-Invert has optimal value lesser than zero (allowing for a numerical tolerance) at the current candidate value $\bar{y}$. (ia) If so, build the corresponding unified cut and add it as a lazy constraint to RMP. (ib) Otherwise, $\bar{y}$ is admissible for IP and $f^T\bar{y} + \widehat{\textit{V}}(\bar{y})$ is a valid upper bound for the optimal value of IP. Update the incumbent solution with the candidate value if appropriate. Step (ii) is identical to that of Section~\ref{sec:Benders algo}.

\subsection{Unified cut of Seo et al.}
\label{sec:Seo unified cut}

\cite{SeoEtAl2022} establish their proposition of a Benders cut on these first principles: (i) Among all optimality cuts (\ref{MP:2}) find the one intersecting the segment between the current RMP solution $(\bar{y}, \bar{\eta})$ and a reference point $(\tilde{y}, \tilde{\eta})$ located in $epi(\widehat{\textit{V}})$ at the point that is closest to $epi(\widehat{\textit{V}})$. (ii) Among all feasibility cuts (\ref{MP:3}), find the one that satisfies the same condition as in (i). (iii) Between the cuts selected in (i) and (ii), select the cut whose intersecting point is closest to $epi(\widehat{\textit{V}})$. \cite{SeoEtAl2022} demonstrate that the cut with the desired properties can be found with this CGP:
\paragraph{CGP-S:}
\begin{align}    
  \max_{\hat{\pi}, \hat{\pi}_0}\; & (G\tilde{y} - q)^T \hat{\pi} - \tilde{\eta} \hat{\pi}_0 \label{CGP-S:1}\\
  \text{subject to } & \hat{\pi}^T G\, d^y - d^{\eta} \hat{\pi}_0 = 1, \label{CGP-S:2}\\
  & -\hat{\pi}^T R - c^T \hat{\pi}_0 \leq 0, \label{CGP-S:3}\\
  & \hat{\pi}, \hat{\pi}_0 \geq 0, \label{CGP-S:4}
  \end{align}
where $(d^y, d^{\eta}) :\equiv (\bar{y} - \tilde{y}, \bar{\eta} - \tilde{\eta})$ and the reference point $(\tilde{y}, \tilde{\eta})$ must satisfy all Benders cuts. (Hence, $(\tilde{y}, \tilde{\eta})$ may simply solve a continuous relaxation of IP.) \label{eq:8.11} If CGP-S is infeasible, no cut is required at $(\bar{y}, \bar{\eta})$. Furthermore, if the optimal value of CGP-S is smaller or equal to -1, no cut is required since $(\bar{y}, \bar{\eta}) \in epi(\widehat{\textit{V}})$. The resulting cut is formally identical to \ref{CGP-FBC_cut:1}. If it is required, it is guaranteed to be supporting $epi(\widehat{\textit{V}})$, and if the solution $(\hat{\pi}, \hat{\pi}_0)$ is unique, it is also supporting a facet of $epi(\widehat{\textit{V}})$. If $(\tilde{y}, \tilde{\eta}) \in relint\, conv(epi(\widehat{\textit{V}}))$ the cut is Pareto-optimal.

\paragraph{Benders iteration.}
(i) Verify that CGP-S is feasible and that its optimal value is lesser than one (allowing for a numerical tolerance) at the current candidate value $\bar{y}$. Steps (ia), (ib) and (ii) are identical to those of Section~\ref{sec:Benders algo}.

\subsection{Unified cut of Hosseini-Turner}
\label{sec:Hosseini-Turner unified cut}

The perspective slightly differs from that of Sections~\ref{sec:Benders algo},~\ref{sec:unified cuts of FBC},~\ref{sec:Seo unified cut} and additional definitions are required. Defining IP's value function as $\textit{V}(y) :\equiv \min_{x \geq 0} \lbrace f^Ty + c^T x \mid Rx \leq q - Gy \rbrace$ and $\textit{V}$'s epigraph as
\begin{equation}
epi(\textit{V}) = \lbrace(y, \theta) \in conv(\mathcal{Y}) \times \mathbb{R} \mid \exists \,x \geq 0 :  Rx \leq q - Gy, f^Ty + c^T x \leq \theta\rbrace,
\end{equation}
IP can be expressed as $\min \lbrace \theta \mid (y, \theta) \in epi(\textit{V}) \cap (\mathcal{Y} \times \mathbb{R})  \rbrace.$
By Farkas' lemma, $(\bar{y}, \bar{\theta}) \in epi(\textit{V})$ iff the following certificate-generating problem has a bounded optimal value (equal to zero):
\begin{align}
\max_{\pi, \pi_0}\; & \pi^T(G\bar{y} - q) + \pi_0 (f^T \bar{y} -\bar{\theta}) \\
\text{subject to }\; & (\pi, \pi_0) \in \Pi :\equiv \lbrace \pi \geq 0, \pi_0 \geq 0 \mid \pi^T R + \pi_0 c^T \geq 0 \rbrace.
\end{align}
If $(\bar{y}, \bar{\theta}) \notin epi(\textit{V}), \exists \, (\pi, \pi_0) \in \Pi \mid \pi^T(G\bar{y} - q) + \pi_0 (f^T \bar{y} -\bar{\theta}) > 0$ whence $(\bar{y}, \bar{\theta})$ breaches the valid constraint $\pi^T(Gy - q) + \pi_0 (f^T y -\theta) \leq 0$. 
Hence, $epi(\textit{V}) = \lbrace (y, \theta) : (y, \theta) \in \mathcal{H}(\pi,\pi_0), \forall \, (\pi, \pi_0) \in \Pi \rbrace$, where $\mathcal{H}(\pi, \pi_0) :\equiv \lbrace (y, \theta) : \pi^T(Gy -q) + \pi_0 (f^T y -\theta) \leq 0 \rbrace$, and a new MP with unified cuts may be defined as
\begin{align}
\min_{y, \theta}\; & \theta \\
\text{subject to }\; & (y, \theta) \in \mathcal{H}(\pi, \pi_0), \forall (\pi, \pi_0) \in \Pi, \\
 & y \in \mathcal{Y}.
\end{align}

The so-called \emph{deepest cuts} proposed by \cite{HossTurn2021, HossTurn2025} progressively outer-approximate $epi(\textit{V})$. For each current solution $(\bar{y}, \bar{\theta}) \in epi(\textit{V})$ yielded by RMP, the cut is defined by the half-space $\mathcal{H}(\pi, \pi_0)$ supporting $epi(\textit{V})$ whose distance to $(\bar{y}, \bar{\theta})$ is largest according to a particular $l_q$ norm. This distance is given by 
\begin{equation} \label{hoss-turn:17}
\delta_{l_q}(\bar{y}, \bar{\theta}\mid \pi, \pi_0 \in \Pi) :\equiv \frac{\pi^T(G\bar{y} - q) + \pi_0(f^T\bar{y} - \bar{\theta})}{\lVert \pi_0f^T + \pi^TG, \pi_0\rVert_p},
\end{equation}
where $\frac{1}{p} + \frac{1}{q} = 1$
and the selected cut is the solution to 
$\delta_{l_q}^*(\bar{y}, \bar{\theta}) = \max_{(\pi, \pi_0) \in \Pi} \delta_{l_q}(\bar{y}, \bar{\theta}\mid \pi, \pi_0).$ Furthermore, $\delta_{l_q}^*(\bar{y}, \bar{\theta})$ is also the $l_q$ distance of $(\bar{y}, \bar{\theta})$ to $epi(\textit{V})$: 
$\delta_{l_q}^*(\bar{y}, \bar{\theta}) = \min_{(y, \theta) \in epi(\mathcal{V})}  \lVert (y - \bar{y}, \theta - \bar{\theta}) \rVert_q.$ In other words, the solution to the RHS problem is the $l_q$ projection of $(\bar{y}, \bar{\theta})$ onto $epi(\textit{V})$. This solution is the point where the cut supports $epi(\textit{V})$.

The CGP is as follows:
\paragraph{CGP-H:}
 \begin{align}
  \max_{\pi, \pi_0}\; & \pi^T(-q + G\bar{y}) + \pi_0(f^T \bar{y} - \bar{\theta}) \label{CGP-H:1}\\ 
  \text{subject to } & \|\pi_0 f^T + \pi^T G,\; \pi_0\|_p \leq 1, \label{CGP-H:2}\\
  & -\pi^T R - \pi_0\, c^T \leq 0, \label{CGP-H:3}\\
  & \pi, \pi_0 \geq 0 \label{CGP-H:4}
\end{align}
and the resulting cut is
\begin{equation}
(-q + Gy)^T \pi + (f^T y - \theta) \pi_0 \leq 0.
\end{equation}
We consider three particular cases, assuming that $G$ has $I$ rows and $J$ columns:
\begin{itemize}
\item When $p=\infty$, i.e. $q = 1$, the constraint $\lVert \pi_0f^T + \pi^TG, \pi_0\rVert_{\infty} \leq 1$ can be replaced with $2J$ linear constraints $-1 \leq \pi_0 f_j + \pi^T G_{.j} \leq 1, \, \forall \, j = 1\dots J$ and $\pi_0 \leq1$. 

\item When $p = 1$, i.e. $q = \infty$, the constraint $\lVert \pi_0f^T + \pi^TG, \pi_0\rVert_1 = \pi_0 + \sum_{j=1}^J \vert \pi_0f_j + \pi^T G_{.j} \vert \leq 1$ can be replaced with an additional variable $\tau \in \mathbb{R}_+^{J}$ and $2J$ constraints $- \tau \leq \pi_0 f + \pi^TG \leq \tau$.

\item One can obtain a relaxation of $l_p, p = 1$, by applying the triangle inequality to $\lVert \pi_0f^T + \pi^TG, \pi_0\rVert_1 = \pi_0 + \sum_{j=1}^J \vert \pi_0f_j + \pi^T G_{.j} \vert =  \pi_0 + \sum_{j=1}^J \vert \pi_0f_j + \sum_{i = 1}^I \pi_i G_{ij} \vert$:
\begin{equation}
\lVert \pi_0 f^T + \pi^TG, \pi_0 \rVert_1 
\leq \pi_0 (1 + \sum_{j=1}^J \vert f_j
\vert) + \sum_{i =1}^I \pi_i \sum_{j=1}^J \vert G_{ij} \vert.
\end{equation}
\end{itemize}

\paragraph{Benders iteration.}
(i) Verify that CGP-H is feasible at the current candidate value $\bar{y}$. (ia) If so, build the corresponding unified cut and add it as a lazy constraint to RMP. (ib) Otherwise, $\bar{y}$ is admissible for IP and $\textit{V}(\bar{y})$ is a valid upper bound for the optimal value of IP. Update the incumbent solution with the candidate value if appropriate. Step (ii) is identical to that of Section~\ref{sec:Benders algo}.

\subsection{Unified cut of Glomb et al.}
\label{sec:Glomb unified cut}

\cite{GlombEtAl2026} introduce the concept of $\beta$-\emph{dominance} between Benders cuts that applies to both feasibility and optimality cuts: Let $\beta$ be an upper bound for the value of IP (a usual choice is the current incumbent value yielded by the branch-and-Benders-cut algorithm). Define the \emph{solution candidate set} associated with a cut $\omega$ as the set of points in $conv(\mathcal{Y})$ that are feasible in RMP and lead to values that are not larger than $\beta$. Cut $\omega_1$ dominates cut $\omega_2$ if the solution candidate set associated with $\omega_1$ is strictly included in the solution candidate associated with $\omega_2$. $\beta$-dominance among optimality cuts is implied by dominance between optimality cuts according to Magnanti-Wong dominance. Hence, non-dominated optimality cuts according to $\beta$-dominance will be non-dominated according to Magnanti-Wong dominance. However the concept of $\beta$-dominance extends to both optimality and feasibility cuts.

\cite{GlombEtAl2026} define the \emph{optimal line-shifting} (OLS) cut. Essentially, the OLS cut excludes the longest segment possible, starting from the current solution $\bar{y}$, on a line between $\bar{y}$ and a reference point $\tilde{y}$ in $conv(\mathcal{Y})$. An OLS cut is demonstrated to be non-dominated according to $\beta$-dominance if the end-point of this segment happens to be in $conv(\mathcal{Y})$. However, there is no guarantee that the latter will be satisfied for a particular cut.
OLS cuts are computed with the following CGP:
\paragraph{CGP-G:}
\begin{align}
  \max_{\hat{\pi}, \hat{\pi}_0}\; & \hat{\pi}^T(G\bar{y} - q) + \hat{\pi}_0(f^T \bar{y} - \beta) \label{CGP-G:1}\\
  \text{subject to }\; & \hat{\pi}^T G(\bar{y} - \tilde{y}) + \hat{\pi}_0 f^T(\bar{y} - \tilde{y}) = 1, \label{CGP-G:2}\\
  & -\hat{\pi}^T R \leq \hat{\pi}_0 c^T, \label{CGP-G:3}\\
  & \hat{\pi} \geq 0, \hat{\pi}_0 \geq 0. \label{CGP-G:4}
\end{align}
The resulting OLS cut is $\hat{\pi}^T(G y - q) \leq \hat{\pi}_0 \eta$. Notice that CGP-G is similar to the adaptive version of CGP-FBC-Invert with the difference that  $\beta-f^Ty$ replaces $\eta$.

Two prerequisites impinge on the availability of OLS cuts: \emph{First}, the current solution $\bar{y}$ must not have led to an improvement of the incumbent value $\beta$. \emph{Second}, a non-trivial regularity condition must also be satisfied. That these prerequisites must be ascertained at every Benders iteration entails that the algorithm governing the application of OLS cuts is of higher complexity than those examined thus far in connection with other unified cuts. In addition to CGP-G, this algorithm involves computations with DSP as well as this modification of DSP where constraint (\ref{DSP-G:4}) is added:

\paragraph{DSP-G:}
\begin{align}
\max_{\hat{\pi}}\; & (G\bar{y} - q)^T \hat{\pi} \\
\text{subject to }\; & -R^T \hat{\pi} \leq c, \\
& (G\tilde{y} - q)^T \hat{\pi} = \beta - f^T \tilde{y}, \label{DSP-G:4}\\
 & \hat{\pi} \geq 0.
\end{align}

\paragraph{Benders iteration.} See \cite{GlombEtAl2026}, Algorithm~1.

\section{The multi-commodity capacitated fixed cost network design problem}
\label{sec:MCFNDP}

Network design problems \citep{CraiEtAl2021} occur in critical, fast-growing areas such as supply chains, transportation, and telecommunications and pursuing the development and assessment of methods for the solution of the increasingly complex NP-hard problems arising in these areas is equally critical. In this endeavor, the linear deterministic MCFNDP naturally plays a central role as a testing instrument by offering a stylized but near-realistic setting. Interest for the linear deterministic MCFNDP as a testing instrument also stems from two other motivations: (i) If the variety of the MCFNDP instances in the testing bench is sufficiently wide and deep, experimental outcomes about the comparative performances of a selection of competing methods can be viewed as indicative of the relative performances of these methods in connection with other comparable high-complexity MILP architectures. (ii) Results obtained with a testing bench made up of instances of a deterministic MCFNDP can be viewed as indicative of the relative performances that would be achieved in a stochastic setting where the stochastic elements have finite support. Hence, once the leading performers have been identified with a deterministic MCFNDP, their performance may be examined closely in the stochastic setting. This option is important in view of the large and often prohibitive resources required by stochastic experiments of realistic scales.

The next two sections respectively present the canonical element-explicit, arc-based statement of the MCFNDP and its matrix representation as a specialization of IP.

\subsection{Element-explicit statement of the MCFNDP}
\label{sec:element-specific statement MCFNDP}

This section is adapted from  \cite{LarsEtAl2023}. Let $\mathcal{G}=(\mathcal{N},\mathcal{A})$ denote a graph composed of arcs $(i,j)\in \mathcal{A}$ and nodes $i,j\in\mathcal{N}$. With each arc $(i,j)\in \mathcal{A}$ are associated a fixed cost $f_{ij}$ and a capacity $u_{ij}\geq 0$ limiting the total amount of flow on the arc. Demands are associated with commodities $k \in \mathcal{K}$ and defined over the nodes of the graph. The unit costs for using any arc $(i,j) \in \mathcal{A}$ are denoted by $c_{ij}^k$. It is commonly assumed that each commodity is characterized by a single origin node $O(k)\in \mathcal{N}$, a single destination node $D(k)\in\mathcal{N}$ and a volume $d^k \geq 0$.  In this setting, the net outgoing flow of commodity $k$ at node $i$ is given by
\[
w_i^k =
\begin{cases}
d^k, & \text{if } i = O(k), \\
-d^k, & \text{if } i = D(k), \\
0, & \text{otherwise}.
\end{cases}
\]
The problem may also include commodity-specific capacities $b_{ij}^k\geq 0$ limiting the flow of specific commodities on the arcs.

The linear deterministic MCFNDP comprises two sets of decision variables: binary design variables $y_{ij},~(i,j)\in\mathcal{A}$, and continuous multi-commodity flow variables $x_{ij}^k\geq 0,~(i,j)\in\mathcal{A},k\in\mathcal{K}$. The arc-based formulation of the MCFNDP is as follows:
\begin{align}
\min_{\boldsymbol{y},\boldsymbol{x}}\; & \sum_{(i,j) \in \mathcal{A}} f_{ij} y_{ij} + \sum_{k\in \mathcal{K}} \sum_{(i,j) \in \mathcal{A}} c^k_{ij} x^k_{ij}  \label{MCFNDP:1}\\ 
\text{subject to }\; & \sum_{j \in \mathcal{N}_i^+} x_{ij}^k - \sum_{j \in \mathcal{N}_i^-} x_{ji}^k = w_i^k, & \forall i \in \mathcal{N}, \forall k \in \mathcal{K}, \label{MCFNDP:2}\\
  & \sum_{k \in \mathcal{K}} x_{ij}^k \leq u_{ij} y_{ij}, & \forall (i,j) \in \mathcal{A}, \label{MCFNDP:3}\\
  & x_{ij}^k \leq b_{ij}^k y_{ij}, & \forall (i,j) \in \mathcal{A}, \forall k \in \mathcal{K}, \label{MCFNDP:4}\\
 & x_{ij}^k \geq 0, & \forall (i,j) \in \mathcal{A}, \forall k \in \mathcal{K}, \label{MCFNDP:5}\\
& y_{ij} \in \lbrace 0,1 \rbrace, & \forall (i,j) \in \mathcal{A}, \label{MMCFNDP:6}
\end{align} 
where $\boldsymbol{y},\boldsymbol{x}$ respectively stand for all design variables $y_{ij},~(i,j)\in\mathcal{A}$, and all flow variables $x_{ij}^k,~(i,j)\in\mathcal{A},k\in\mathcal{K}$. The objective function (\ref{MCFNDP:1}) minimizes the total costs expressed as a summation of fixed arc opening costs and a summation of flow-dependent transportation costs. Constraints~(\ref{MCFNDP:2}) enforce flow conservation at each node $i\in \mathcal{N}$, where $\mathcal{N}_i^+$ and $\mathcal{N}_i^-$ identify the successor and predecessor nodes of $i$, respectively.
Constraints~(\ref{MCFNDP:3}) enforce overall capacity limits and act as linking constraints. Finally,~(\ref{MCFNDP:4}) are commodity-specific capacity constraints (a.k.a. strong capacity constraints) restricting the flow of commodity $k\in\mathcal{K}$ on arc $(i,j)$. Their inclusion does not change the optimal solution of the MCFNDP but it can improve the strength of its linear programming relaxation.

\subsection{Matrix representation of the MCFNDP}
\label{sec:matrix representation MCFNDP}

The vectors and matrices of the generic IP are specialized as follows to represent the MCFNDP. Let $N :\equiv |\mathcal{N}|$, $A :\equiv |\mathcal{A}|$, $K :\equiv |\mathcal{K}|$. We suppose that arcs $(i,j) \in \mathcal{A}$ are indexed with $a = 1, \ldots, A$.

\begin{itemize}
\item[] $y$: Vector of size $A$ whose elements are indexed with $a = 1, \ldots, A$.
\item[] $x$: Vector of size $AK$ whose elements are indexed by the pairs $(a, k)$, $a = 1, \ldots, A$; $k = 1, \ldots, K$ enumerated in lexicographic order.
\item[] $f$: Vector of same size and indexing as $y$.
\item[] $c$: Vector of same size and indexing as $x$.

\item[]

$G = \begin{bmatrix} G(1A) \\ G(1B) \\ G(2) \\ G(3) \end{bmatrix},
\quad
R = \begin{bmatrix} R(1A) \\ R(1B) \\ R(2) \\ R(3) \end{bmatrix},
\quad
q = \begin{bmatrix} g(1A) \\ g(1B) \\ g(2) \\ g(3) \end{bmatrix},$

whose submatrices are defined as follows:
\item[] $G(1A)$, $G(1B)$: Null matrices of size $NK \times A$.

\item[] $R(1A)$: Matrix of size $NK \times AK$ whose rows are indexed with the pairs $(i, k)$, $i = 1, \ldots, N$; $k = 1, \ldots, K$, enumerated in lexicographic order. In the row corresponding to $(i, k)$, all elements are null except that those whose columns correspond to arcs leaving node $i$ for commodity $k$ are equal to 1 and those whose columns correspond to arcs entering node $i$ for commodity $k$ are equal to -1. 

\item[] $R(1B)$: Equal to $-R(1A)$.

\item[] $q(1A)$: Vector of size $NK$ indexed with the pairs $(i, k)$, $i = 1, \ldots, N$; $k = 1, \ldots, K$ enumerated in lexicographic order, whose element corresponding to $(i, k)$ is equal to $w_i^k$.

\item[] $q(1B)$: Equal to $-q(1A)$.

\item[] $G(2)$: Matrix of size $A \times A$ whose rows and columns are indexed with the ordered arcs. In row $a$ corresponding to arc $a$, all elements are null except that corresponding to column $a$ is equal to $-u_a$.

\item[] $R(2)$: Matrix of size $A \times AK$ whose rows are indexed with the ordered arcs and whose columns are indexed with the pairs $(a, k)$ ordered in lexicographic order. In row corresponding to a, all elements are null except those corresponding to pair $(a, k)$, $k = 1, \ldots, K$ are equal to 1.

\item[] $q(2)$: Null vector of size $A$.

\item[] $G(3)$: Matrix of size $AK \times A$ whose rows are indexed with $(a, k)$, $a = 1, \ldots, A$; $k = 1, \ldots, K$, enumerated in lexicographic order and whose columns are indexed with the ordered arcs. In row corresponding to pair $(a, k)$, all elements are null except that corresponding to column $a$ is equal to $-b_a^k$.

\item[] $R(3)$: Identity matrix of size $AK$.

\item[] $q(3)$: Null vector of size $AK$.
\end{itemize}

If all strong capacity constraints in (\ref{MCFNDP:4}) are omitted, blocks $G(3)$, $R(3)$ and $q(3)$ are excluded from $G$, $R$ and $q$. If some or all strong capacity constraints are included in (\ref{MCFNDP:4}), the corresponding rows of $G(3)$, $R(3)$, $g(3)$ are included in $G$, $R$ and $q$.

\paragraph{Element-explicit expansions.} The Appendix expresses the CGPs and the cuts of Sections~\ref{sec:Benders algo}, \ref{sec:unified cuts} in terms of the individual elements of their constituent matrices and vectors when the latter are instantiated according to the present section.

\section{Description of experiments}
\label{sec:description of experiments}

In total, we evaluate and compare the performances of 52 computation methods (50 Benders + 2 direct). Each method is characterized by an algorithmic family and a set of hyperparameters shaping its exact behavior. 

\paragraph{Algorithmic families.}
For each algorithmic family, we state its short-hand identifier used in Section~\ref{sec:experimental results} when reporting experimental results, include a brief description and point to the section describing the relevant CGP and Benders iteration, if applicable.

\begin{enumerate}
\item[] ILP\_flow: Direct branch-and-cut solution of IP performed by a solver.
\item[] Std: Standard Benders decomposition with distinct feasibility and optimality cuts based on DSP (Section~\ref{sec:Benders algo}).
\item[] Std\_Decomp: Uses the built-in flow decomposition implemented in SMS++ to solve the subproblem of the MCFNDP with a MILP solver. In our context, this amounts to standard Benders decomposition with distinct feasibility and optimality cuts based on PSP (Section~\ref{sec:Benders algo}).
\item[] MW: Application of Magnanti-Wong optimality cut \citep{MagnWong1981} along with standard feasibility and optimality cuts (Section~\ref{sec:Magnanti-Wong optim cut}).  
\item[] Papa: Application of Papadakos optimality cut \citep{Papa2008} along with standard feasibility and optimality cuts (Section~\ref{sec:Papadakos optim cut}).
\item[] Fisch: Application of unified cut according to \cite{FiscEtAl2010} (Section~\ref{sec:unified cuts of FBC}).
\item[] StatBS: Application of unified static cut according to \cite{BranStur2021} (Section~\ref{sec:unified cuts of FBC}).
\item[] AdaptBS: Application of unified adaptive cut according to \cite{BranStur2021} (Section~\ref{sec:unified cuts of FBC}).
\item[] Seo: Application of unified cut according to \cite{SeoEtAl2022} (Section~\ref{sec:Seo unified cut}).
\item[] HT: Application of unified cut according to \cite{HossTurn2021, HossTurn2025} (Section~\ref{sec:Hosseini-Turner unified cut}).
\item[] Glomb: Application of unified cut according to \cite{GlombEtAl2026} (Section~\ref{sec:Glomb unified cut}).
\end{enumerate}

\paragraph{Hyperparameters.}
For each hyperparameter, we state its short-hand identifier and include a brief description indicating purpose, admissible values and relevant algorithmic family(ies).

\begin{enumerate}
\item[] F: Indicates inclusion of strong capacity constraints in MCFNDP. Relevant to all algorithmic families.
\item[] U: Indicates use of CGP-FBC rather than CGP-FBC-Invert. Relevant to Fisch and StatBS.
\item[] G: Value of relative gap that must be attained in preliminary application of Fisch ([0, 1]). Relevant to AdaptBS and Glomb. 
\item[] M: Number of successive admissible integer solution(s) that must be identified in preliminary application of Fisch (strictly positive integer). Relevant to AdaptBS and Glomb. (Recall that methods belonging to these algorithmic families require an incumbent solution to start.)
\item[] E: Method of computation for reference point (2: solution of continuously relaxed IP calculated only once, initially; 4: solution of continuously relaxed IP, augmented with cuts previously and calculated at each Benders iteration). Relevant to Seo.
\item[] R: Underlying distance function (0: $l_\infty$; 1: $l_1$; 2: upper bound on $l_1$). Relevant to HT.
\item[] P: Common value of all coordinates of initial reference point ([0, 1]). Relevant to MW and Papa.
\item[] L: Mixing weight updating the reference point by combining it with current solution ([0, 1]). Relevant to MW and Papa.
\item[] H: Mixing weight updating the reference point by combining it with current solution ([0, 1]). Relevant to Glomb.
\item[] T: Relative gap threshold ([0, 1]). The computation process is interrupted and the incumbent solution is considered optimal when the relative gap attains this value. This is set to 0.0001 for all algorithmic families and hyperparameter sets. Hence, mention of this hyperparameter can be omitted for brevity in all method identifiers without risk of confusion.
\item[] A: Also add a standard optimality cut whenever a Magnanti-Wong optimality cut is added. Relevant to MW.
\end{enumerate}

\paragraph{MCFNDP instances.}
As a testing bench for evaluating the computational characteristics of all methods, we select the Canad R set of MCFNDP instances \citep{Canad}. The Canad R set comprises 162 generated problems featuring systematically varying numbers of nodes, arcs and commodities and systematically varying fixed costs and capacities, thus promoting diversified structural characteristics and creating a broad, progressive spectrum of computational complexities. Computation times with the fastest methods range from within seconds to the time limit whereas the relative gap is still far from zero. The Canad R set has been used since its inception as a testing bench in several empirical analyses by several groups of authors (see, e.g., \cite{ChouEtAl2016,CraiEtAl2001,HewiEtAl2010}). Arcs of the Canad R instances can be interpreted from the source files as (origin, destination) or (destination, origin), thus giving rise to two distinct sets of feasible MCFNDP instances with cardinalities respectively equal to 162 and 153. We follow the custom of reading destinations first. The 153 valid Canad R instances are hence denoted $r\,x.y$, where $x = 01,\dots,03$, $y=1,\dots,6$ and $x = 04,\dots,18$, $y=1,\dots,9$.

\paragraph{Computational apparatus.}
Computing programs are written in the C++ language supported by the open-source GCC 15.1.0 toolchain. The C++ open-source SMS++ 0.4.0 structured mathematical modeling and optimization library \citep{SMS++} performs the backend computations. We run two parallel sets of experiments where the supporting MILP solver called by SMS++ is respectively Gurobi 12.0.0 \citep{Guro2026} and CPLEX 22.1.2 \citep{CPLE2026}. Computations are performed on a workstation equipped with two AMD EPYC 7763 processors of 64 cores each and 2TB of ECC DDR4-3200MHz RAM, managed by the Linux RHEL operating system. Each experiment is allocated a single core and a time limit of 10h.

\section{Measuring performance} 
\label{sec:measuring performance}

The next sections describe in turn (i) the raw measurements that are collected throughout the experiments, (ii) the bespoke evaluation methods based on relative differences and rankings of raw measurements that are used to conduct a first round of evaluation of all methods, (iii) the performance profiles \citep{DolaMore2002} used in a second round of evaluation to discriminate between leaders identified in the first round.

\subsection{Raw measurements}
\label{sec:raw measurements}

The following measurements are collected in relation to the application of each one of the computation methods included in the analysis to each one of the 153 Canad R MCFNDP instances: status at exit, total computation time, value of lowest upper bound, value of highest lower bound, relative gap, total number of nodes generated, total number of CGPs solved, total numbers of feasibility and optimality cuts generated, total numbers of feasible CGPs with positive, negative and null optimal values, total number of infeasible CGPs, total numbers of positively and negatively unbounded CGPs, total sum of times spent in CGPs, total sum of squared times spent in CGPs, maximum resident memory size, maximum virtual memory size.

\subsection{Aggregating performance over all problem instances}
\label{sec:aggregating performance}

Comparing the overall performance of the computational methods over all problem instances (PIs) must account for two caveats:

\paragraph{First caveat: Availability and relevance of raw measurements depend on outcomes.}
Terminal values for the performance measurements listed in Section~\ref{sec:raw measurements} are not available and meaningful for every computational outcome. On the one hand, when computations are interrupted by the time limit, only the terminal values of the relative gap and the relative differences to the best upper and lower bounds are genuinely meaningful. On the other hand, terminal gap and relative differences to the best upper and lower bounds are trivially negligible when computation completes successfully before or at time limit. This must be accounted for when (i) summarizing the performance attained by the individual methods over the PIs and (ii) comparing summaries between methods.

We respond as follows: First, divide the set of Canad PIs between Blocks A and B: PIs in A are solved by at least one of the Benders methods under examination before or at time limit. The remainder of the PIs are in B. Second, divide the computational (method, PI) pairs between Blocks A1 and A2: Pairs in A1 are made up of a method solving the paired PI in Block A before or at time limit. Pairs in A2 are made up of a method that does not solve the paired PI in Block A before or at time limit. We say that a (method, PI) pair is in Block B when the paired PI belongs to Block B. Equipped with the partition of all (method, PI) pairs between Blocks A1, A2 and B, we can meaningfully define the aggregate performance of each method over the PIs in a block. Notice that this partition of the (method, PI) pairs between Blocks A1, A2 and B is tied to the particular set of Benders methods under examination and to the particular setting chosen for the experiments, namely the set of PIs making up the testing bench and the computational apparatus. 

\paragraph{Second caveat: Raw measurements are related to size and complexity of PIs.} Clearly, the terminal value of each performance measurement listed in Section~\ref{sec:raw measurements} is related to the size and complexity of the particular PI about which they are collected. Since the experiments involve PIs of widely varying sizes and complexities, this must be accounted for when (i) summarizing the performance attained by the individual method over the PIs and (ii) comparing the summaries between methods.

We respond by avoiding comparisons based directly on the raw performance measurements and by instead comparing methods based on their relative differences or their rankings with respect to the performance measurements.

\subsubsection{Comparing performances based on relative differences}
\label{sec:compare based on rel diffs}
Let us see now how relative differences can alternatively be used. Suppose that computational methods and PIs are respectively indexed with $m = 1,\dots,M$ and $p=1,\dots,P$ and that $\mu_{mp}$ is one of the measurements enumerated in Section~\ref{sec:raw measurements} taken in connection with the pair $(m,p)$. The relative difference is defined as $\delta_{mp} :\equiv (\mu_{mp} - \min_m \mu_{mp}) / \min_m \mu_{mp}$ if small values of $\mu_{mp}$ are desirable, and as $\delta_{mp} :\equiv (\max_m \mu_{mp} -\mu_{mp}) / \max_m \mu_{mp}$ otherwise. Again, let us focus as an example on Block A1 and computation time. For each particular value of PI appearing in at least one (method, PI) pair in Block A1, find the smallest computation time over all pairs where the same value of PI also appears. Then, calculate the relative difference to this figure for each pair. Repeat for every particular value of PI. We end up with distributions of relative differences for each particular method encountered among the pairs in Block A1. We can then calculate standard summary statistics such as mean and quantiles for each one of the distributions associated with a particular method. Comparison of the aggregate performance of the methods over the PIs can then be based on these statistics. Clearly, beyond this example similar calculations can be performed in regard to any particular performance measure relevant to any particular Block A1, A2 or B.

\subsubsection{Comparing performances based on rankings}
\label{sec:compare based on rankings}
Let us see how rankings can be used to compare for example the (method, PI) pairs in Block A1 with respect to computation time. Thus, in regard to all (method, PI) pairs in Block A1 sharing a common PI, say Canad $r\,a.b$, we rank the methods according to computation time in increasing order. Then we repeat in this fashion for each particular value of PI encountered among the pairs in Block A1. We end up with distributions of rankings for each particular method encountered among the pairs in Block A1. Calculation of summary statistics for each one of the distributions associated with a particular method and use of these statistics to compare their performance can proceed in connection with the use of the rankings as described in Section~\ref{sec:compare based on rel diffs}. Here also, similar calculations can be performed beyond this example in regard to any particular performance measure relevant to any particular Block A1, A2 or B.

\subsection{Performance profiles}
Performance profiles \citep{DolaMore2002} offer an illuminating graphical representation when a small number of lists of scores must be compared.  
Let the performance ratio be defined as $r_{mp} :\equiv \frac{\mu_{mp}}{\min_{m}\mu_{mp}}$, where $m$ and $p$ stand for a method and a PI, $\mu_{mp}$ 
is the value of raw measurement $\mu$ resulting from the application of method $m$ to problem $p$, and smaller values of $\mu$ are desirable. The performance profile of a particular method $m$ in regard to raw measurement $\mu$ is defined 
as $\rho_m(\tau) = \frac{1}{n_p} card\lbrace p \in P \mid r_{mp} \leq \tau\rbrace$, where $\mathcal{P}$ is the set of all PIs, of cardinality $n_p$. Thus, the performance profile associated with method $m$ indicates the proportion of all PIs for which $m$ would offer the best score if the score of the best method were artificially multiplied by $\tau$. Thus, for every particular measurement taken in relation to the (method, PI) pairs occurring in one of the Blocks A1, A2 and B, there would be up to 52 lists of up to 153 scores to handle. Performance profiles could not be used in these circumstances since including a large number of profiles on a single plot (up to 52 here) would make them illegible. However, they can be compared over a small number of leading methods that have been identified from examination of relative differences and/or rankings as in Sections~\ref{sec:compare based on rel diffs} and \ref{sec:compare based on rankings}.

\section{Experimental results}
\label{sec:experimental results}

This section reports on the experimental results. Two parallel groups of experiments are conducted. In the first one, SMS++ calls Gurobi to run the computations. This group of experiments includes all 52 methods under consideration. The results are examined in Section~\ref{sec:report on computations with Gurobi}. In the second group of experiments, SMS++ calls CPLEX. The results are examined in Section~\ref{sec:report on computations with CPLEX}. Table~\ref{tab:all experiments run} below summarizes all experiments run with Gurobi and with CPLEX.

We present selected excerpts of the complete set of output data files. The latter is available in the repository located at \textbf{https://github.com/larseeri/unified-Benders-cuts-compared}. The included \textbf{README.md} file supplies the requisite information about its exact content.


\begin{table}[htbp]
\centering
\small
\begin{tblr}{
  colspec = {l l c c c c},
  cell{1}{1,2} = {r=2}{l},
  cell{1}{3} = {c=2}{c},
  cell{1}{5} = {c=2}{c},
  row{1,2} = {font=\bfseries},
}
\hline
Algorithmic family & Hyperparameters & Gurobi & & CPLEX & \\
       &               & $-$F & $+$F & $-$F & $+$F \\
\hline
Direct    & ILP\_flow                 & \ev   & \rep  & \ev & \rep \\
\hline
Std       & Std                       & \ev   & \ev   &     &      \\
          & Std\_Decomp               & \ev   & \rep  &     & \rep \\
\hline
Fisch     & Fisch                     & \ev   & \ev  & \ev  & \ev \\
          & U\_Fisch                  & \ev   & \rep  &  \ev  & \rep \\
\hline
StatBS    & StatBS                    & \ev   & \lead & \ev   & \lead \\
          & U\_StatBS                 & \ev   & \lead &  \ev   & \lead \\
\hline
AdaptBS   & G-0.1\_M-1                & \ev   & \ev   &     &      \\
          & G-0.5\_M-1                & \ev   & \ev   &     &      \\
          & M-1                       & \ev   & \ev  &     & \rep \\
\hline
Seo       & E-2                       & \ev   & \rep  &     & \rep \\
          & E-4                       & \ev   & \ev   &     &      \\
\hline
HT        & R-0                       & \ev   & \ev   &     &      \\
          & R-1                       & \ev   & \lead &  \ev   & \lead \\
          & R-2                       & \ev   & \lead &  \ev   & \lead \\
\hline
MW        & P-1.0\_L-0.5\_A           & \ev   & \rep  &     & \rep \\
\hline
Papa      & P-1.0\_L-0.5              & \ev   & \ev   &     &      \\
\hline
Glomb     & H-0.1\_G-0.1\_M-1         & \ev   & \ev   &     &      \\
          & H-0.1\_G-0.25\_M-1        & \ev   & \ev   &     &      \\
          & H-0.1\_G-0.75\_M-1        & \ev   & \ev   &     &      \\
          & H-0.5\_G-0.1\_M-1         & \rep  & \ev   & \rep &      \\
          & H-0.5\_G-0.25\_M-1        & \ev   & \ev   &     &      \\
          & H-0.5\_G-0.75\_M-1        & \ev   & \ev   &     &      \\
          & H-0.9\_G-0.1\_M-1         & \ev   & \ev   &     &      \\
          & H-0.9\_G-0.25\_M-1        & \ev   & \ev   &     &      \\
          & H-0.9\_G-0.75\_M-1        & \ev   & \ev   &     &      \\
\hline
\end{tblr}
\caption{All experiments (26 base methods, each without ($-$F)/with ($+$F) strong
capacity constraints, run on 153 Canad R instances). \ev: evaluated; \rep: evaluated and reported as a family
representative in Tables~\ref{tab:rel distances blocks A1 A2 Gurobi}--\ref{tab:rel distances block B CPLEX};
\lead: found a leading method and reported in Tables~\ref{tab:rel distances blocks A1 A2 Gurobi}--\ref{tab:rel distances block B CPLEX}
and Figures~\ref{fig:mosaic performance profiles block A1 Gurobi}--\ref{fig:mosaic performance profiles block B Gurobi}. All 52 methods are evaluated under Gurobi; a representative subset is evaluated under CPLEX.}
\label{tab:all experiments run}
\end{table}

\subsection{Main set of experiments: Gurobi}
\label{sec:report on computations with Gurobi}

\subsubsection{Methods investigated}

Using Gurobi as the supporting solver of the SMS++ library, we run each one of the 26 pairs of methods whose identifiers are listed below on each one of the 153 Canad R instances of MCFNDP.
\begin{itemize}
\item[] (ILP\_flow, F\_ILP\_flow) 
\item[] (Std, F\_Std)
\item[] (Std\_Decomp, F\_Std\_Decomp)
\item[] (Fisch, F\_Fisch)
\item[] (U\_Fisch, F\_U\_Fisch)
\item[] (StatBS, F\_StatBS)
\item[] (U\_StatBS, F\_U\_StatBS)
\item[] (G-0.1\_M-1\_AdaptBS, F\_G-0.1\_M-1\_AdaptBS)
\item[] (G-0.5\_M-1\_AdaptBS, F\_G-0.5\_M-1\_AdaptBS)
\item[] (M-1\_AdaptBS, F\_M-1\_AdaptBS)
\item[] (E-2\_Seo, F\_E-2\_Seo)
\item[] (E-4\_Seo, F\_E-4\_Seo)
\item[] (R-0\_HT, F\_R-0\_HT)
\item[] (R-1\_HT, F\_R-1\_HT)
\item[] (R-2\_HT, F\_R-2\_HT)
\item[] (P-1.0\_L-0.5\_A\_MW, F\_P-1.0\_L-0.5\_A\_MW)
\item[] (P-1.0\_L-0.5\_Papa, F\_P-1.0\_L-0.5\_Papa)
\item[] (H-0.1\_G-0.1\_M-1\_Glomb, F\_H-0.1\_G-0.1\_M-1\_Glomb)
\item[]	(H-0.1\_G-0.25\_M-1\_Glomb, F\_H-0.1\_G-0.25\_M-1\_Glomb)
\item[] (H-0.1\_G-0.75\_M-1\_Glomb, F\_H-0.1\_G-0.75\_M-1\_Glomb)
\item[] (H-0.5\_G-0.1\_M-1\_Glomb, F\_H-0.5\_G-0.1\_M-1\_Glomb)
\item[] (H-0.5\_G-0.25\_M-1\_Glomb, F\_H-0.5\_G-0.25\_M-1\_Glomb)
\item[] (H-0.5\_G-0.75\_M-1\_Glomb, F\_H-0.5\_G-0.75\_M-1\_Glomb)
\item[] (H-0.9\_G-0.1\_M-1\_Glomb, F\_H-0.9\_G-0.1\_M-1\_Glomb)
\item[] (H-0.9\_G-0.25\_M-1\_Glomb, F\_H-0.9\_G-0.25\_M-1\_Glomb)
\item[] (H-0.9\_G-0.75\_M-1\_Glomb, F\_H-0.9\_G-0.75\_M-1\_Glomb)
\end{itemize}
Although the identifiers are quite transparent about the exact nature of the methods they represent, some clarifications are necessary or at least useful: (i) Methods in a pair differ only with respect to inclusion of the strong capacity constraints (parameter F). (ii) Presence of parameter G in connection with the Glomb et al. algorithmic family indicates that the specified relative gap must be attained through a preliminary application of Fisch before control is turned over to Algorithm~1 of \cite{GlombEtAl2026}. (Notice that parameter M ensures that the requisite incumbent value is made available.) We implement three values for G: 0.1, 0.25 and 0.75. Value 0.1 is meant to imitate the \emph{hybrid strategy} proposed in \cite{GlombEtAl2026}. This is the most favorable of the two strategies proposed by the authors. Hence, if the resulting performance is dominated by those of concurrent methods, the ensuing conclusion will be clear. Values 0.25 and 0.75 serve to verify how well the Glomb et al. algorithm can perform on its own earlier in the solution process. (iii) The exact value of parameter H in connection with the Glomb et al. family is left open by \cite{GlombEtAl2026}. We implement values of 0.1, 0.5 and 0.9 to verify if and how they impact performance. (iv) Presence of parameter G in connection with the AdaptBS algorithmic family indicates that the specified relative gap must be attained through a preliminary application of Fisch before control is turned over to the procedure involving adaptive Brandenberg-Stursberg cuts delineated in Section~\ref{sec:unified cuts of FBC}. We consider three cases: no preliminary application of Fisch as in the original description of \cite{BranStur2021}, inclusion of G with values of 0.1 and 0.5. The latter cases are meant to verify if the finer but more time-consuming calculations of AdaptBS may be more useful in the later stage of the solution process. 

\subsubsection{Results}

\paragraph{Algorithmic settings and computational characteristics.}
Computations are performed with Gurobi's parameters standing at their default values. Selection of both root and node algorithms is therefore set at automatic mode (\emph{dual simplex} since a single core is available). We find out that increasing the value of \emph{NumericFocus} does not remedy the large number of invalid computations observed for some of the methods (usually characterized by a failure to find a valid integer solution within the maximum time limit). We also find out that achieving successful computations with the methods belonging to the adaptBS, Seo et al. and Glomb et al. algorithmic families hinges on delicate adjustments to the magnitude and signs of the numerical tolerances involved in responding to the outcome of the relevant CGPs. These tolerances serve to slightly shift the theoretical bounds triggering the addition of cuts. We notice that these algorithmic families present either CGPs that include computed differences between optimization outcomes (adaptBS, Seo) or a complex algorithm with several steps and CGPs (Glomb). This might have the effect of compounding the sources of numerical errors. In the end, the numerical sensitivity of adaptBS, Seo and Glomb complicates their implementation and may suggest a lack of robustness in further applications.   

\paragraph{Relative differences and rankings.}
In accordance with the method of Section~\ref{sec:compare based on rel diffs}, 
Tables~\ref{tab:rel distances blocks A1 A2 Gurobi} and \ref{tab:rel distances block B Gurobi} contain excerpts from the complete results about the relative differences and report on a selection of methods stemming from these preliminary findings:
(i) Inclusion of strong capacity constraints leads to important reductions of computation times everywhere except for the methods belonging to the Glomb et al. family. (ii) The direct method F\_ILP\_flow is the fastest overall. (iii) Methods belonging to the Glomb et al. family fail or perform poorly when strong capacity constraints are included. (iv) Method F\_Std is dominated by F\_Std\_Decomp. (v) F\_G-0.1\_M-1\_AdaptBS, F\_G-0.5\_M-1\_AdaptBS and F\_M-1\_AdaptBS fail to complete their computations before the time limit for several PIs in Block A. We select F\_M-1\_AdaptBS to illustrate the typical performance achieved by methods belonging to the AdaptBS family. (vi) In similarity with methods in the AdaptBS family, F\_E-2\_Seo and F\_E-4\_Seo fail to complete their computations before the time limit for several PIs in Block A. We select F\_E-2\_Seo to illustrate the typical performance achieved by methods belonging to the Seo family. (vii) F\_R-1\_HT and F\_R-2\_HT by far dominate F\_R-0\_HT. (viii) F\_P-1.0\_L-0.5\_A\_MW dominates F\_P-1.0\_L-0.5\_Papa. Still, F\_P-1.0\_L-0.5\_A\_MW fails to complete computations before the time limit for several PIs in Block A. (ix) Values of parameters H and G in connection with the methods belonging to the Glomb et al. family have marginal impacts on performance, which is poor. We select H-0.5\_G-0.1\_M-1\_Glomb to illustrate the typical performance achieved by the methods belonging to the Glomb et al. family.

Summing up, examination of relative differences in Tables~\ref{tab:rel distances blocks A1 A2 Gurobi} and \ref{tab:rel distances block B Gurobi} indicates that among the Benders methods, the following stand out: F\_U\_StatBS, F\_StatBS, F\_R-1\_HT, F\_R-2\_HT. They stand out in Block A in regard to the following considerations: number of valid computations finishing before time limit; computation time; numbers of nodes, subproblems and cuts; time spent in subproblems. They also stand out in Block B in regard to number of valid computations, relative gap and distance to best upper and lower bounds. F\_U\_Fisch, F\_Fisch, F\_Std\_Decomp stand far behind. Analysis of rankings according to the method of Section~\ref{sec:compare based on rankings} leads to similar conclusions. 

\begin{table}[htbp]
\centering
\normalsize

\begin{adjustbox}{width=\textwidth,totalheight=\textheight, keepaspectratio}

\begin{tblr}{colspec={||l|rrrrrrrrrrr||}}

\hline
\SetCell[r=2]{c}  method        & F\_ &F\_ &F\_U\_ &F\_ & F\_U\_ &F\_M-1\_ &F\_E-2\_ &F\_R-1\_ &F\_R-2\_ &F\_P-1.0\_L-0.5\_A\_ &H-0.5\_G-0.1\_M-1\_ \\
& ILP\_flow & Std\_Decomp &  Fisch &  StatBS & StatBS & AdaptBS & Seo & HT & HT & MW & Glomb \\   
\hline
\SetCell[c=12]{c} Block A: Report on MCFNDPs solved within time limit (36000s) through at least one Benders method &       &       &       &       &       &       &       &       &       &   &     \\
\hline     
   numb valid        & 90 & 89    & 90    & 90    & 90 & 89    & 88    & 89    & 90    & 89    & 89 \\
\hline
\SetCell[c=12]{c}  Block A1: Computations with methods finishing within time limit (36000s) about MCFNDPs in Block A &       &       &       &       &       &       &       &       &       &       &    \\
\hline
   numb valid       & 90 & 79    & 81    & 89    & 90 & 74    & 76    & 88    & 89    & 74    & 73 \\
\hline
\SetCell[c=12]{l}  Relative distance to minimum computation time &       &       &       &       &       &       &       &       &       &       &         \\
          mean  & -0.78 & 12.68 & 18.43 & 2.57  & 2.44 & 26.90 & 36.47 & 2.87  & 2.58  & 36.79 & 29.82 \\
          std err & 0.03 & 2.63  & 5.77  & 0.68  & 0.87 & 10.70 & 13.65 & 0.62  & 0.66  & 8.75  & 9.25 \\
           0.25 q & -0.98 & 0.33  & 0.87  & 0.32  & 0.00 & 2.36  & 1.32  & 0.45  & 0.28  & 4.45  & 1.35 \\
           0.5 q &-0.92 & 3.54  & 1.81  & 0.72  & 0.32 & 4.71  & 3.66  & 1.04  & 0.71  & 14.63 & 3.21 \\
           0.75 q & -0.73 & 15.17 & 9.31  & 1.48  & 0.77 & 11.27 & 13.31 & 1.96  & 1.58  & 28.75 & 9.10 \\
\SetCell[c=12]{l}  Relative distance to minimum number of nodes &       &       &       &       &       &       &       &       &       &       &         \\
           mean  & -0.52 & 65.83 & 35.03 & 2.66  & 2.49 & 52.20 & 4728.88 & 2.80  & 3.06  & 87.03 & 71.17 \\
           std err & 0.19 & 19.42 & 12.98 & 0.89  & 0.97 & 17.44 & 4335.54 & 1.10  & 1.15  & 27.09 & 20.20 \\
           0.25 q & -1.00 & 2.86  & 0.76  & 0.09  & 0.02 & 2.04  & 0.00  & 0.00  & 0.00  & 3.45  & 1.96 \\
           0.5 q & -0.99 & 7.59  & 2.32  & 0.49  & 0.35 & 4.50  & 2.18  & 0.26  & 0.32  & 10.60 & 6.79 \\
           0.75 q & -0.84 & 24.89 & 10.38 & 1.35  & 1.04 & 17.27 & 10.97 & 0.94  & 1.21  & 37.40 & 29.80 \\
\SetCell[c=12]{l}   Relative distance to minimum number of subproblems &       &       &       &       &       &       &       &       &       &       &    \\
           mean     &  none  & 3.55  & 3.00  & 0.19  & 0.23 & 2.91  & 3.79  & 0.23  & 0.25  & 5.23  & 5.15 \\
           std err    &  none  & 0.44  & 0.54  & 0.03  & 0.03 & 0.33  & 0.71  & 0.03  & 0.03  & 0.72  & 0.98 \\
           0.25 q    &  none  & 1.26  & 0.61  & 0.00  & 0.00 & 0.99  & 0.66  & 0.01  & 0.03  & 1.84  & 0.66 \\
           0.5 q     &  none  & 2.40  & 1.18  & 0.10  & 0.16 & 2.29  & 1.34  & 0.15  & 0.19  & 3.65  & 2.10 \\
           0.75 q    &  none  & 4.24  & 3.36  & 0.29  & 0.34 & 3.99  & 4.21  & 0.34  & 0.34  & 5.99  & 4.83 \\
\SetCell[c=12]{l}   Relative distance to minimum number of feasibility cuts &       &       &       &       &       &       &       &       &       &       &     \\
           mean    &  none  & undef & undef & undef & undef  & undef & undef & undef & undef & undef & undef \\
           std err    &  none  & undef & undef & undef & undef   & undef & undef & undef & undef & undef & undef \\
           0.25 q   &  none  & 7.64  & 1.33  & 0.05  & 0.09 & 1.50  & 0.00  & 0.00  & 0.00  & 6.84  & 1.59 \\
           0.5 q     &  none  & 17.83 & 2.44  & 0.27  & 0.31 & 2.42  & 0.76  & 0.24  & 0.22  & 14.30 & 2.88 \\
          0.75 q   &  none  & 35.13 & 4.60  & 0.62  & 0.65 & 4.15  & 2.00  & 0.75  & 0.67  & 31.58 & 4.80 \\
\SetCell[c=12]{l}   Relative distance to minimum number of optimality cuts &       &       &       &       &       &       &       &       &       &       &    \\
          mean    &  none  & 1.56  & 3.76  & 0.54  & 0.54 & 3.84  & 7.30  & 0.59  & 0.62  & 3.81  & 6.56 \\
           std err   &  none  & 0.20  & 0.64  & 0.11  & 0.11 & 0.48  & 1.79  & 0.11  & 0.11  & 0.43  & 1.17 \\
          0.25 q    &  none  & 0.41  & 0.62  & 0.03  & 0.02 & 1.35  & 0.81  & 0.05  & 0.07  & 1.34  & 0.84 \\
           0.5 q     &  none  & 1.07  & 1.49  & 0.17  & 0.24 & 2.90  & 2.19  & 0.29  & 0.28  & 2.86  & 2.87 \\
          0.75 q    &  none  & 2.04  & 3.85  & 0.57  & 0.60 & 4.82  & 6.24  & 0.62  & 0.68  & 5.01  & 7.17 \\
 \SetCell[c=12]{l}   Relative distance to minimum number of cuts &       &       &       &       &       &       &       &       &       &       &     \\
          mean     &  none  & 3.65  & 3.10  & 0.19  & 0.21 & 3.00  & 4.01  & 0.23  & 0.25  & 5.30  & 5.40 \\
          std err    &  none  & 0.47  & 0.56  & 0.03  & 0.03 & 0.34  & 0.73  & 0.03  & 0.03  & 0.74  & 1.02 \\
          0.25 q     & none & 1.27  & 0.62  & 0.00  & 0.00 & 0.94  & 0.72  & 0.03  & 0.03  & 1.88  & 0.73 \\
          0.5 q     &  none & 2.41  & 1.28  & 0.08  & 0.11 & 2.33  & 1.49  & 0.14  & 0.19  & 3.61  & 2.39 \\
          0.75 q    &  none  & 4.28  & 3.39  & 0.26  & 0.35 & 4.09  & 4.71  & 0.31  & 0.38  & 6.06  & 4.91 \\
\SetCell[c=12]{l}   Relative distance to minimum time spent in subproblems &       &       &       &       &       &       &       &       &       &       &    \\
          mean    &  none  & 1.55  & 10.16 & 3.31  & 2.54 & 12.89 & 30.51 & 3.76  & 3.31  & 17.55 & 10.05 \\
          std err    &  none  & 0.43  & 3.80  & 0.70  & 0.88 & 4.35  & 12.53 & 0.62  & 0.68  & 4.22  & 2.99 \\
          0.25 q    &  none  & 0.06  & 0.94  & 0.54  & 0.00 & 2.44  & 1.41  & 1.19  & 0.71  & 3.73  & 0.98 \\
          0.5 q      &  none  & 0.62  & 1.79  & 1.09  & 0.27 & 4.11  & 4.04  & 1.98  & 1.19  & 7.09  & 2.92 \\
          0.75 q    &  none  & 1.70  & 5.56  & 2.70  & 0.92 & 6.75  & 13.21 & 3.52  & 2.55  & 13.07 & 7.61 \\
\hline
 \SetCell[c=12]{c}   Block A2: Computations with methods reaching time limit (36000s) about MCFNDPs in Block A &       &       &        &       &       &       &       &       &       &       &  \\
\hline
   numb valid        & 0 & 10 & 9  & 1  & 0 & 15 & 12 & 1  & 1  & 15 & 16 \\
\hline
\SetCell[c=12]{l}   Relative distance to minimum relative gap &       &       &       &       &       &       &       &       &       &       &       \\
          mean  &  none  & 6.33  & 5.73  & 0.42  & none   & 13.02 & 4.60  & 0.00  & 1.99  & 8.89  & 5.94 \\
          std err &  none  & 2.55  & 2.22  & 0.00  & none   & 7.07  & 1.97  & 0.00  & 0.00  & 3.97  & 1.72 \\
          0.25 q &  none  & 0.42  & 0.59  & 0.42  & none   & 0.45  & 0.00  & 0.00  & 1.99  & 0.41  & 1.13 \\
           0.5 q &  none  & 2.87  & 2.11  & 0.42  & none   & 2.16  & 1.05  & 0.00  & 1.99  & 2.61  & 2.71 \\
           0.75 q &  none  & 10.16 & 10.03 & 0.42  & none   & 10.75 & 6.30  & 0.00  & 1.99  & 7.04  & 9.49 \\
\SetCell[c=12]{l}   Relative distance to minimum upper bound &       &       &       &       &       &       &       &       &       &       &       \\
           mean  &  none  & 0.01  & 0.02  & 0.00  & none   & 0.02  & 0.02  & 0.00  & 0.00  & 0.02  & 0.03 \\
           std err &  none  & 0.01  & 0.01  & 0.00  & none   & 0.01  & 0.01  & 0.00  & 0.00  & 0.01  & 0.01 \\
           0.25 q &  none  & 0.00  & 0.00  & 0.00  & none   & 0.00  & 0.00  & 0.00  & 0.00  & 0.00  & 0.00 \\
           0.5 q &  none  & 0.00  & 0.00  & 0.00  & none   & 0.00  & 0.01  & 0.00  & 0.00  & 0.00  & 0.00 \\
           0.75 q &  none  & 0.01  & 0.04  & 0.00  & none   & 0.02  & 0.03  & 0.00  & 0.00  & 0.02  & 0.06 \\
\SetCell[c=12]{l}   Relative distance to maximum lower bound &       &       &       &       &       &       &       &       &       &       &       \\
           mean  &  none  & 0.07  & 0.05  & 0.00  & none   & 0.03  & 0.01  & 0.00  & 0.02  & 0.09  & 0.05 \\
           std err &  none  & 0.03  & 0.01  & 0.00  & none   & 0.01  & 0.00  & 0.00  & 0.00  & 0.05  & 0.01 \\
           0.25 q &  none  & 0.01  & 0.01  & 0.00  & none   & 0.01  & 0.00  & 0.00  & 0.02  & 0.01  & 0.01 \\
           0.5 q &  none  & 0.04  & 0.03  & 0.00  & none   & 0.03  & 0.01  & 0.00  & 0.02  & 0.01  & 0.03 \\
           0.75 q &  none  & 0.16  & 0.09  & 0.00  & none   & 0.04  & 0.02  & 0.00  & 0.02  & 0.07  & 0.08 \\
\hline
\end{tblr}
\end{adjustbox}
\caption{Relative distances yielded by computations in Blocks A1, A2 with selected methods using Gurobi. The minima and the maximum are calculated over the Benders methods. (undef: undefined since minimum is equal to 0 for some instances.)}
\label{tab:rel distances blocks A1 A2 Gurobi}
\end{table}

\begin{table}[htb]
\centering
\normalsize

\begin{adjustbox}{width=\textwidth,totalheight=\textheight, keepaspectratio}

\begin{tblr}{colspec={||l|rrrrrrrrrrr||}}

\hline
\SetCell[r=2]{l}  method   & F\_ &F\_ &F\_U\_ &F\_ & F\_U\_ &F\_M-1\_ &F\_E-2\_ &F\_R-1\_ &F\_R-2\_ &F\_P-1.0\_L-0.5\_A\_ &H-0.5\_G-0.1\_M-1\_ \\
 & ILP\_flow & Std\_Decomp &  Fisch &  StatBS & StatBS & AdaptBS & Seo & HT & HT & MW & Glomb \\   
\hline
\SetCell[c=12]{c}   Block B: Report on MCFNDPs that CANNOT be solved within time limit (36000s) through at least one Benders method &       &       &        &       &       &       &       &       &       &       &  \\
\hline
numb valid   & 63    & 61	& 63	& 60	& 63     & 61    & 51	& 59	& 61	& 62	& 63 \\
\hline
\SetCell[c=12]{l}   Relative distance to minimum relative gap &       &       &       &       &       &       &       &       &       &       &       \\
          mean  & -0.984  & 4.83	& 5.63	& 1.61	& 0.60	& 2.96	& 1.11	& 1.29	& 1.34	& 15.12	& 6.73
 \\
          std err & 0.004   & 1.09  & 1.99  & 0.76  & 0.27  & 0.63  & 0.28  & 0.42  & 0.55  & 6.58  & 2.25
 \\
          0.25 q 	& -1.00	&  1.32	& 0.68	& 0.16	& 0.00	& 0.85	& 0.05	&  0.29	& 0.16	& 1.69	& 1.09
 \\
           0.5 q    & -0.999  &  1.84	& 1.21	& 0.42	& 0.00	& 1.47	& 0.53	& 0.53	& 0.45  & 3.05	& 1.58
\\
           0.75 q	& -0.983  &  2.94	& 2.65	&  0.86  & 0.31	& 2.98	& 1.20 	& 0.96 	& 0.73    & 7.78  & 3.38
\\
\SetCell[c=12]{l}   Relative distance to minimum upper bound &       &       &       &       &       &       &       &       &       &       &       \\
           mean 	& -0.051   & 0.08	& 0.09	& 0.08	& 0.02	& 0.14	& 0.11	& 0.08	& 0.08	& 0.16	& 0.12
 \\
           std err &  0.006	&0.01		& 0.01	& 0.01	& 0.00	& 0.02	&  0.02	& 0.01	& 0.01	& 0.01	& 0.01
 \\
           0.25 q	&  -0.081	& 0.03	& 0.04	& 0.01	& 0.00	& 0.04	& 0.01	& 0.01	& 0.01	& 0.08	& 0.06
 \\
           0.5 q 	&  -0.039	& 0.07	& 0.07	& 0.06	& 0.00	& 0.13	& 0.09	& 0.05	& 0.04	& 0.14	&0.13
\\
           0.75 q 	&  -0.008 	& 0.13	& 0.13	& 0.14	& 0.03	& 0.20	& 0.16	& 0.17	& 0.13	& 0.22	& 0.19
 \\
\SetCell[c=12]{l}   Relative distance to maximum lower bound &       &       &       &       &       &       &       &       &       &       &       \\
           mean  	& -0.063  & 0.22  & 0.14  & 0.04  & 0.03	& 0.12	& 0.02	& 0.04	& 0.04	& 0.44	& 0.17
\\
           std err & 0.006	& 0.02	& 0.01	& 0.00	& 0.00	& 0.01	& 0.01	& 0.00	& 0.00	& 0.04	& 0.01
 \\
           0.25 q 	& -0.086 & 0.13	& 0.07	& 0.01	& 0.00	& 0.08	& 0.00	& 0.01	& 0.01	& 0.16	& 0.10
 \\
           0.5 q 	& -0.055  & 0.17  & 0.13	& 0.02	& 0.00	& 0.11	& 0.00	& 0.02	& 0.02	& 0.44	& 0.16
 \\
           0.75 q 	&  -0.027    & 0.29	& 0.20	& 0.06	& 0.04	& 0.17	& 0.01	& 0.06	& 0.06	& 0.70	& 0.23
 \\
\hline
\end{tblr}
\end{adjustbox}
\caption{Relative distances yielded by computations in Block B with selected methods using Gurobi. The minima and the maximum are calculated over the Benders methods.}
\label{tab:rel distances block B Gurobi}
\end{table}

\begin{figure}[htb]
    \centering
    \begin{subfigure}[b]{0.32\textwidth}
        \centering
        \includegraphics[width=\linewidth]{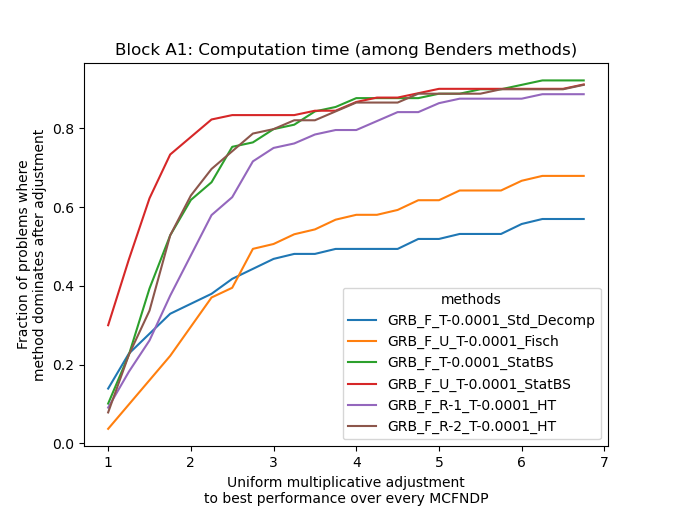}
        \label{fig:mosaic A1_comp_time}
    \end{subfigure}
    \hfill
    \begin{subfigure}[b]{0.32\textwidth}
        \centering
        \includegraphics[width=\linewidth]{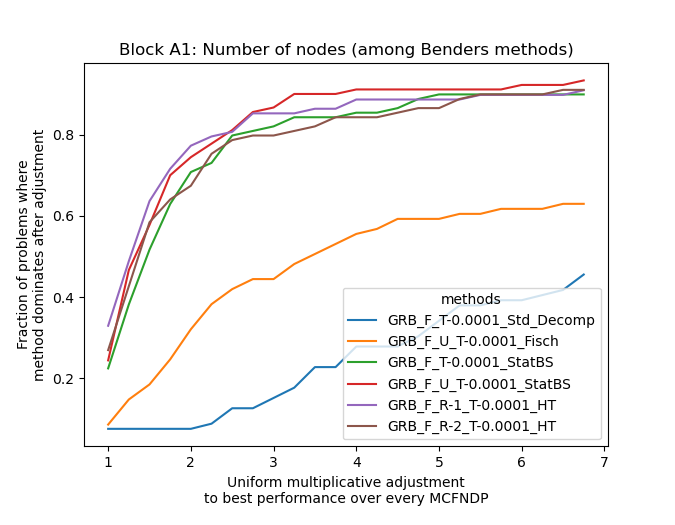}
        \label{fig:mosaic A1_numb_nodes}
    \end{subfigure}
	\hfill
    \begin{subfigure}[b]{0.32\textwidth}
        \centering
        \includegraphics[width=\linewidth]{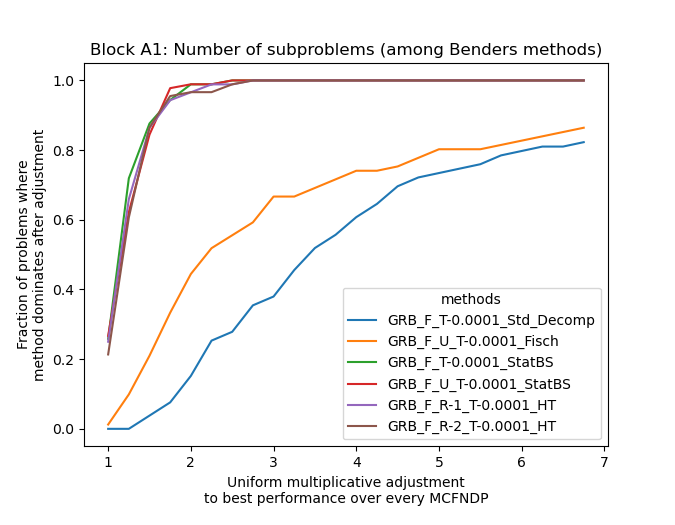}
        \label{fig:mosaic A1_numb_subp}
    \end{subfigure}    
    \begin{subfigure}[b]{0.32\textwidth}
        \centering
        \includegraphics[width=\linewidth]{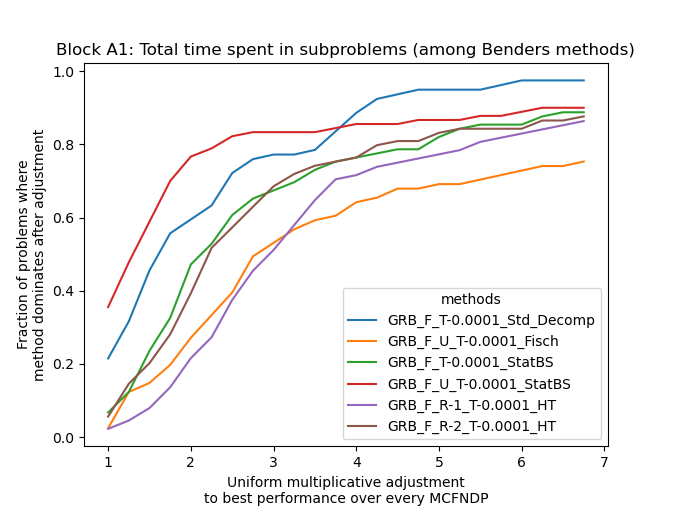}
        \label{fig:mosaic A1_time_subp}
    \end{subfigure}
    \hfill
    \begin{subfigure}[b]{0.32\textwidth}
        \centering
        \includegraphics[width=\linewidth]{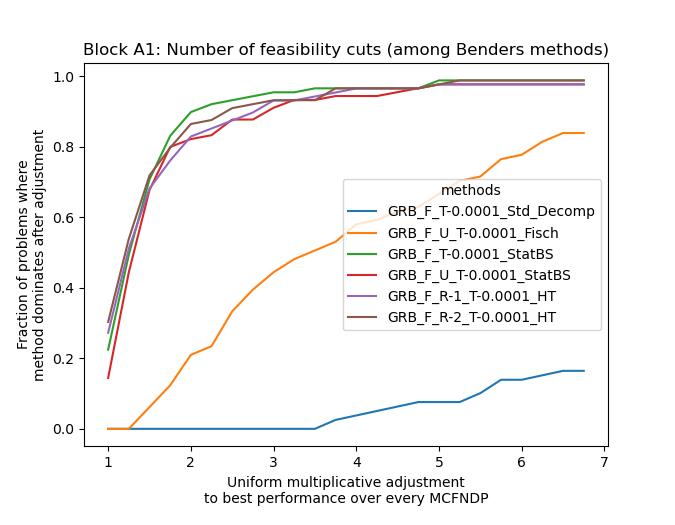}
        \label{fig:mosaic A1_numb_cuts}
    \end{subfigure}
	\hfill
    \begin{subfigure}[b]{0.32\textwidth}
        \centering
        \includegraphics[width=\linewidth]{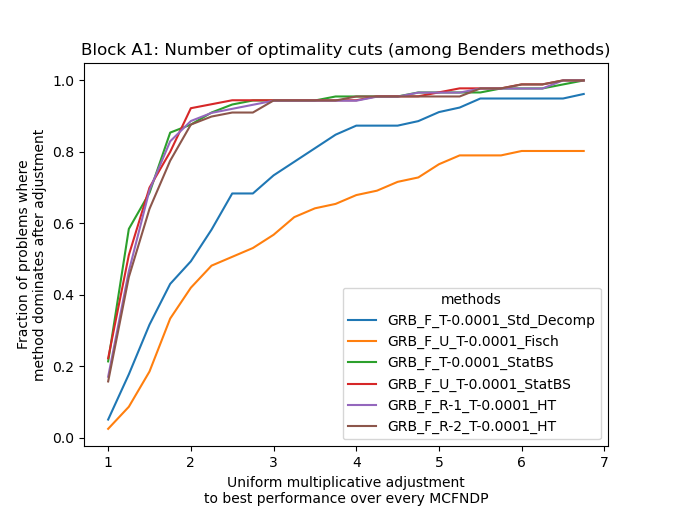}
        \label{fig:mosaic A1_numb_opt_cuts}
    \end{subfigure}	
    \caption{Performance profiles in Block A1 of leading methods}
    \label{fig:mosaic performance profiles block A1 Gurobi}
\end{figure}

\begin{figure}[htb]
    \centering
    \begin{subfigure}[b]{0.32\textwidth}
        \centering
        \includegraphics[width=\linewidth]{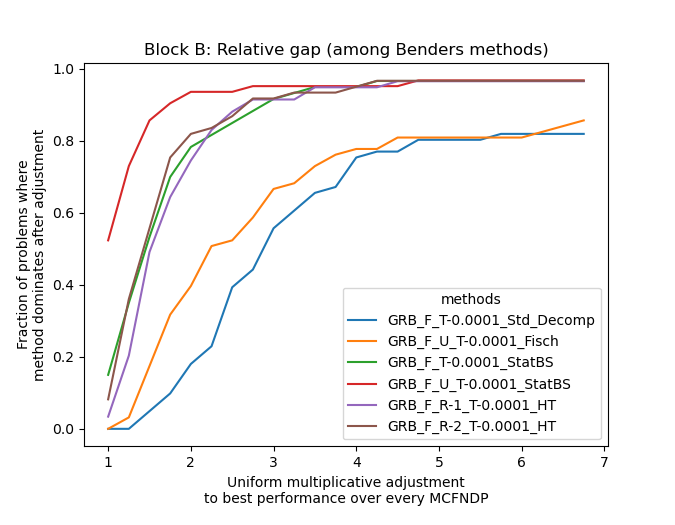}
        \label{fig:mosaic B_rel_gap}
    \end{subfigure}
    \hfill
    \begin{subfigure}[b]{0.32\textwidth}
        \centering
        \includegraphics[width=\linewidth]{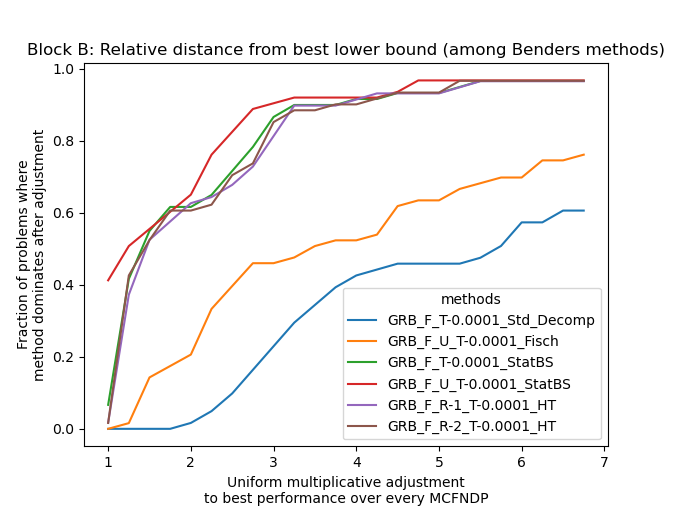}
        \label{fig:mosaic B_rel_dist_from_lwrBd}
    \end{subfigure}
	\hfill
    \begin{subfigure}[b]{0.32\textwidth}
        \centering
        \includegraphics[width=\linewidth]{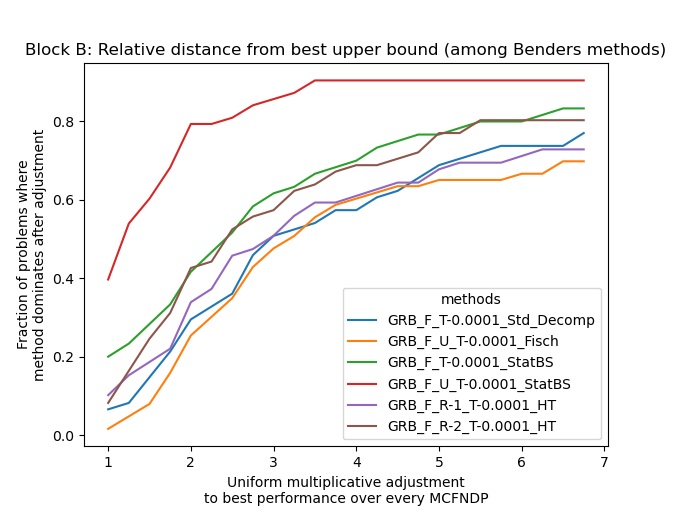}
        \label{fig:mosaic B_rel_dist_from_uprBd}
    \end{subfigure}    
    \caption{Performance profiles in Block B of leading methods}
    \label{fig:mosaic performance profiles block B Gurobi}
\end{figure}

\paragraph{Performance profiles.}
Figures~\ref{fig:mosaic performance profiles block A1 Gurobi} and \ref{fig:mosaic performance profiles block B Gurobi} present performance profiles for the leading methods F\_U\_StatBS, F\_StatBS, F\_R-1\_HT, F\_R-2\_HT and also F\_Std\_Decomp, F\_U\_Fisch in regard to the following measurements: (i) Within Block~A1: computation time, number of nodes, number of subproblems, total time spent in subproblems, number of feasibility and optimality cuts. (ii) Within Block~B: relative gap, relative distances to highest lower bound and lowest upper bound available (calculated with F\_ILP\_flow). Notice that since the leading methods naturally appear in none or very few of the pairs in Block~A2, the related performance profiles are either undefined or present little reliability. Hence, they are omitted. Clearly, F\_U\_StatBS, first, and F\_StatBS, F\_R-2\_HT, F\_R-1\_HT, second, dominate F\_Std\_Decomp and F\_U\_Fisch with respect to computation time in Block~A and relative gap in Block~B. Domination is also manifest in Block~A with respect to numbers of nodes, subproblems and cuts and in Block~B with respect to relative distance from best upper and lower bounds.

\subsection{Second set of experiments: CPLEX}
\label{sec:report on computations with CPLEX}

In the objective of comparing the computational characteristics of Gurobi and CPLEX and to verify their agreement on a set of leading methods, we perform computations with CPLEX as the supporting solver on a subset of the methods examined in Section~\ref{sec:report on computations with Gurobi}. 

\subsubsection{Methods investigated}

We include the leaders F\_U\_StatBS, F\_StatBS, F\_R-1\_HT, F\_R-2\_HT identified in Section~\ref{sec:report on computations with Gurobi} with Gurobi as the supporting solver. To establish a global portrait and enable stronger conclusions we also include the following methods identified in Section~\ref{sec:report on computations with Gurobi}: (i) the direct methods ILP\_flow and F\_ILP\_flow, (ii) the less successful Benders methods F\_U\_Fisch, F\_Fisch, U\_Fisch, Fisch, U\_StatBS, StatBS, R-1\_HT, R-2\_HT, F\_Std\_Decomp, (iii) representatives from  the less successful algorithmic families F\_E-2\_Seo, F\_P-1.0\_L-0.5\_A\_MW, H-0.5\_G-0.1\_M-1\_Glomb, F\_M-1\_AdaptBS. Each one of the aforementioned methods is run on each one of the 153 Canad R MCFNDP instances with CPLEX as the supporting solver.

\subsubsection{Results}

\paragraph{Algorithmic settings and computational characteristics.}
Computations are performed with CPLEX's parameters standing at their default values, with one exception described below. Selection of both root and node algorithms is set at automatic mode (dual simplex since a single core is available). We found out that CPLEX consumes prohibitively large amounts of RAM for node storage under the default setting (compressed nodes stored in RAM). These amounts very often exceed by two orders of magnitude those required by Gurobi, even for median-sized instances (i.e., less than 5GB on average vs up to 300GB for computation involving a single method and a single instance). We therefore set CPLEX's parameter CPXPARAM\_MIP\_Strategy\_File to 2: nodes are stored on disk uncompressed.

\paragraph{Relative differences and rankings.}

Recall from Section~\ref{sec:aggregating performance} that the partition between Blocks A1, A2 and B is tied to the set of Benders methods under examination, to the set of PIs making up the testing bench and to the computational apparatus. Hence, Blocks A1, A2 and B of this section are not identical to those of Section~\ref{sec:report on computations with Gurobi}. 
Examination of relative differences, excerpted in Tables~\ref{tab:rel distances blocks A1 A2 CPLEX} and \ref{tab:rel distances block B CPLEX}, yields the following observations echoing those of Section~\ref{sec:report on computations with Gurobi}: (i) Inclusion of strong capacity constraints leads to major improvements in performance. (ii) The direct methods F\_ILP\_flow and ILP\_flow are the fastest overall. (iii) Performances of the Benders methods follow essentially the same order as with Gurobi and F\_U\_StatBS, F\_StatBS, F\_R-1\_HT, F\_R-2\_HT still stand out. (iv) The number of valid solutions returned by F\_ILP\_flow and ILP\_flow are comparable with those observed with Gurobi.

However, Table~\ref{tab:rel distances blocks A1 A2 CPLEX} also illustrates one important difference from Section~\ref{sec:report on computations with Gurobi}: (v) The number of invalid solutions returned by the Benders methods is significantly larger with CPLEX than with Gurobi. Upon closer examination, this is caused by CPLEX yielding suboptimal solutions in Block A1 for several MCFNDP instances.
Since these differing numbers of invalid solutions between Gurobi and CPLEX might blur the comparisons between the relative differences calculated with each solver, we repeat the calculations of Tables~\ref{tab:rel distances blocks A1 A2 CPLEX} and \ref{tab:rel distances block B CPLEX} while ensuring that the statistics computed for each method are based on the same set of valid instances. (The additional results are available in the repository.) We conclude that observations (i) to (iv) above hold. Findings similar to (i)-(v) result from the examination of rankings.

\paragraph{Computation times.}

CPLEX's computation times are significantly longer than Gurobi's. As an illustration, we calculate the weighted average of computation times for CPLEX and for Gurobi over the instances that have valid solutions and also finish computation before the time limit with both CPLEX and Gurobi. The ratios between these figures for the direct methods F\_ILP\_flow, ILP\_flow using branch-and-cut and for the Benders methods F\_U\_StatBS, F\_StatBS, F\_R-1\_HT, F\_R-2\_HT using branch-and-Benders-cut are respectively equal to 1.69, 1.42, 1.78, 1.74, 4.46, 2.52 in favor of Gurobi. (These additional results are calculated from the output files available in the repository.)

\begin{table}[htbp]
\centering
\normalsize

\begin{adjustbox}{width=\textwidth,totalheight=\textheight, keepaspectratio}

\begin{tblr}{colspec={||l|rrrrrrrrrrr||}}

\hline
\SetCell[r=2]{c}  method        & F\_ &F\_ &F\_U\_ &F\_ & F\_U\_ &F\_M-1\_ &F\_E-2\_ &F\_R-1\_ &F\_R-2\_ &F\_P-1.0\_L-0.5\_A\_ &H-0.5\_G-0.1\_M-1\_ \\
& ILP\_flow & Std\_Decomp &  Fisch &  StatBS & StatBS & AdaptBS & Seo & HT & HT & MW & Glomb \\   
\hline
\SetCell[c=12]{c} Block A: Report on MCFNDPs solved within time limit (36000s) through at least one Benders method &       &       &       &       &       &       &       &       &       &   &     \\
\hline     
   numb valid      & 89	&77	&85	&73	&86	&85	&18	&75	&70	&80	&84 \\
\hline
\SetCell[c=12]{c}  Block A1: Computations with methods finishing within time limit (36000s) about MCFNDPs in Block A &       &       &       &       &       &       &       &       &       &       &    \\
\hline
   numb valid       & 89	&64	&73	&71	&86	&67	&6	&70	&67	&60	&70 \\
\hline
\SetCell[c=12]{l}  Relative distance to minimum computation time &       &       &       &       &       &       &       &       &       &       &         \\
          mean  &  -0.84	&40.57	&20.77	&2.83	&2.19	&49.67	&19.57	&8.77	&3.22	&97.53	&29.38 \\
          std err & 0.03	&11.76	&6.46	&0.93	&0.49	&21.98	&13.43	&3.16	&1.09	&36.11	&11.13 \\
           0.25 q &  -0.99	&0.00	&1.46	&0.00	&0.16	&1.56	&3.62	&0.83	&0.04	&3.97	&0.58 \\
           0.5 q & -0.93	&6.32	&4.21	&0.22	&0.50	&4.35	&5.78	&1.73	&0.29	&20.46	&2.97 \\
           0.75 q &  -0.82	&34.37	&12.68	&0.73	&1.61	&19.75	&7.30	&4.07	&1.06	&75.93	&11.08 \\
\SetCell[c=12]{l}  Relative distance to minimum number of nodes &       &       &       &       &       &       &       &       &       &       &         \\
           mean  &  -0.99	&undef	&undef	&0.65	&undef	&undef	&undef	&undef	&0.97	&undef	&undef \\
           std err & 0.00	&undef	&undef	&0.13	&undef	&undef	&undef	&undef	&0.26	&undef	&undef \\
           0.25 q & -1.00	&4.05	&0.53	&0.05	&0.00	&2.09	&3.57	&0.18	&0.01	&3.61	&0.93 \\
           0.5 q &  -1.00	&11.92	&1.14	&0.26	&0.14	&4.38	&17.42	&0.37	&0.33	&13.83	&2.85 \\
           0.75 q & -0.98	&30.45	&7.04	&0.65	&0.58	&13.85	&undef	&0.90	&0.72	&35.67	&13.39 \\
\SetCell[c=12]{l}   Relative distance to minimum number of subproblems &       &       &       &       &       &       &       &       &       &       &    \\
           mean     &  none&  5.32	&2.32	&0.77	&0.08	&5.13	&49.63	&0.75	&1.59	&8.22	&11.35 \\
           std err  &  none  &  0.83	&0.49	&0.13	&0.02	&1.46	&31.01	&0.13	&0.30	&0.87	&2.49 \\
           0.25 q   &  none &  1.61	&0.30	&0.12	&0.00	&1.38	&10.74	&0.13	&0.19	&2.70	&1.25 \\
           0.5 q    &  none &  4.11	&0.70	&0.34	&0.00	&2.41	&20.56	&0.38	&0.71	&6.39	&3.85 \\
           0.75 q   &  none &  6.41	&2.32	&0.99	&0.09	&4.49	&26.44	&1.00	&1.40	&11.81	&9.47 \\
\SetCell[c=12]{l}   Relative distance to minimum number of feasibility cuts &       &       &       &       &       &       &       &       &       &       &     \\
           mean    &  none &undef	&undef	&undef	&undef	&undef	&undef	&undef	&undef	&undef	&undef \\
           std err   &  none &undef	&undef	&undef	&undef	&undef	&undef	&undef	&undef	&undef	&undef \\
           0.25 q   &  none&   5.92	&0.67	&0.00	&0.00	&0.56	&97.90	&0.00	&0.10	&6.38	&1.06 \\
           0.5 q    &  none &  20.75	&1.43	&0.17	&0.12	&1.24	&118.80	&0.18	&0.33	&21.10	&1.94 \\
          0.75 q   &  none&  60.74	&2.40	&0.50	&0.30	&2.15	&undef	&0.53	&0.59	&52.35	&4.25 \\
\SetCell[c=12]{l}   Relative distance to minimum number of optimality cuts &       &       &       &       &       &       &       &       &       &       &    \\
          mean    &  none&  3.86	&2.50	&1.08	&0.14	&6.77	&38.53	&1.03	&2.91	&5.19	&12.99 \\
           std err  &  none &  0.70	&0.52	&0.18	&0.03	&2.29	&25.41	&0.18	&0.84	&0.95	&2.81 \\
          0.25 q    &  none&  0.93	&0.30	&0.15	&0.00	&1.97	&9.61	&0.19	&0.21	&1.15	&1.86 \\
           0.5 q    &  none &  2.73	&0.80	&0.52	&0.00	&2.89	&12.62	&0.47	&0.83	&2.98	&4.39 \\
          0.75 q   &  none &  4.51	&2.32	&1.44	&0.12	&4.85	&17.07	&1.24	&1.52	&6.82	&11.01 \\
 \SetCell[c=12]{l}   Relative distance to minimum number of cuts &       &       &       &       &       &       &       &       &       &       &     \\
          mean     &  none&  5.56	&2.40	&0.81	&0.07	&5.56	&55.05	&0.79	&1.69	&6.21	&11.95 \\
          std err   &  none &  0.88	&0.51	&0.13	&0.02	&1.71	&34.44	&0.14	&0.32	&0.70	&2.61 \\
          0.25 q    &  none & 1.62	&0.34	&0.12	&0.00	&71.44	&11.95	&0.12	&0.19	&1.82	&1.39 \\
          0.5 q    &  none &  4.20	&0.70	&0.35	&0.00	&2.46	&23.34	&0.39	&0.72	&5.46	&4.14 \\
          0.75 q   &  none &  7.00	&2.57	&1.04	&0.06	&4.61	&28.58	&0.96	&1.47	&9.04	&10.39 \\
\SetCell[c=12]{l}   Relative distance to minimum time spent in subproblems &       &       &       &       &       &       &       &       &       &       &    \\
          mean   &  none & 2.40	&8.48	&3.01	&2.58	&22.25	&23.43	&9.67	&3.42	&21.36	&9.08 \\
          std err   &  none &  0.89	&1.38	&0.95	&0.51	&12.17	&17.08	&3.17	&1.10	&3.53	&2.79 \\
          0.25 q    &  none&  0.00	&2.02	&0.00	&0.28	&0.97	&3.34	&1.16	&0.02	&3.83	&0.60 \\
          0.5 q    &  none  &  0.40	&4.11	&0.31	&0.92	&1.79	&5.54	&2.53	&0.35	&10.94	&2.11 \\
          0.75 q    &  none&  2.49	&8.81	&0.98	&2.26	&5.12	&7.31	&5.47	&1.49	&23.48	&6.40 \\
\hline
 \SetCell[c=12]{c}   Block A2: Computations with methods reaching time limit (36000s) about MCFNDPs in Block A &       &       &        &       &       &       &       &       &       &       &  \\
\hline
   numb valid        & 0	&13	&12	&2	&0	&18	&12	&5	&3	&20	&14
 \\
\hline
\SetCell[c=12]{l}   Relative distance to minimum relative gap &       &       &       &       &       &       &       &       &       &       &       \\
          mean  &  none&   31.19	&3.26	&0.32	&none	&14.58	&10.72	&5.16	&0.79	&23.81	&8.95 \\
          std err &  none&  14.90	&2.06	&0.23	&none	&10.24	&7.53	&3.81	&0.53	&10.29	&5.15 \\
          0.25 q &  none&  0.00	&0.21	&0.16	&none	&0.00	&0.00	&0.16	&0.15	&0.20	&0.91 \\
           0.5 q &  none&  4.67	&0.54	&0.32	&none	&2.62	&0.98	&1.54	&0.30	&1.24	&2.17 \\
           0.75 q &  none&  37.98	&1.33	&0.48	&none	&7.53	&4.07	&1.97	&1.19	&20.22	&3.48 \\
\SetCell[c=12]{l}   Relative distance to minimum upper bound &       &       &       &       &       &       &       &       &       &       &       \\
           mean  &  none&  0.10	&0.01	&0.00	&none	&0.01	&0.14	&0.01	&0.00	&0.08	&0.02 \\
           std err &  none&  0.07	&0.00	&0.00	&none	&0.00	&0.07	&0.01	&0.00	&0.05	&0.00 \\
           0.25 q &  none&  0.00	&0.00	&0.00	&none	&0.00	&0.00	&0.00	&0.00	&0.00	&0.00 \\
           0.5 q &  none&  0.01	&0.00	&0.00	&none	&0.00	&0.02	&0.00	&0.00	&0.02	&0.01 \\
           0.75 q &  none&  0.06	&0.02	&0.00	&none	&0.02	&0.10	&0.01	&0.01	&0.04	&0.03 \\
\SetCell[c=12]{l}   Relative distance to maximum lower bound &       &       &       &       &       &       &       &       &       &       &       \\
           mean  &  none&  0.32	&0.05	&0.00	&none	&0.04	&0.10	&0.04	&0.00	&0.20	&0.04 \\
           std err &  none&  0.11	&0.01	&0.00	&none	&0.01	&0.07	&0.03	&0.00	&0.07	&0.01 \\
           0.25 q &  none&  0.01	&0.01	&0.00	&none	&0.00	&0.00	&0.00	&0.00	&0.01	&0.01 \\
           0.5 q &  none&  0.04	&0.04	&0.00	&none	&0.03	&0.00	&0.01	&0.00	&0.03	&0.03 \\
           0.75 q &  none&  0.65	&0.06	&0.00	&none	&0.06	&0.03	&0.01	&0.00	&0.22	&0.07 \\
\hline
\end{tblr}
\end{adjustbox}
\caption{Relative distances yielded by computations in Blocks A1, A2 with selected methods using CPLEX. The minima and the maximum are calculated over the Benders methods. (undef: undefined since minimum is equal to 0 for some instances.)}
\label{tab:rel distances blocks A1 A2 CPLEX}
\end{table}

\begin{table}[htb]
\centering
\normalsize

\begin{adjustbox}{width=\textwidth,totalheight=\textheight, keepaspectratio}

\begin{tblr}{colspec={||l|rrrrrrrrrrr||}}

\hline
\SetCell[r=2]{l}  method   & F\_ &F\_ &F\_U\_ &F\_ & F\_U\_ &F\_M-1\_ &F\_E-2\_ &F\_R-1\_ &F\_R-2\_ &F\_P-1.0\_L-0.5\_A\_ &H-0.5\_G-0.1\_M-1\_ \\
 & ILP\_flow & Std\_Decomp &  Fisch &  StatBS & StatBS & AdaptBS & Seo & HT & HT & MW & Glomb \\   
\hline
\SetCell[c=12]{c}   Block B: Report on MCFNDPs that CANNOT be solved within time limit (36000s) through at least one Benders method &       &       &        &       &       &       &       &       &       &       &  \\
\hline
numb valid   & 64	&63	&61	&60	&60	&48	&24	&53	&60	&62	&64
    \\
\hline
\SetCell[c=12]{l}   Relative distance to minimum relative gap &       &       &       &       &       &       &       &       &       &       &       \\
          mean  &-0.99	&11.55	&3.52	&0.74	&0.38	&7.45	&3.80	&1.49	&1.14	&17.43 &3.81
\\
          std err &0.00	&3.61	&0.70	&0.34	&0.11	&2.45	&2.23	&0.69	&0.73	&4.61 &0.89
\\
          0.25 q 	& -1.00	&1.74	&0.76	&0.02	&0.00	&1.16	&0.00	&0.17	&0.01	&2.22 &0.65
\\
           0.5 q    & -1.00	&4.34	&1.59	&0.11	&0.05	&1.88	&0.30	&0.33	&0.08	&5.87 &1.53
\\
           0.75 q	& -0.98	&8.03	&3.00	&0.51	&0.24	&4.46	&2.02	&0.81	&0.36	&13.16 &3.29
\\
\SetCell[c=12]{l}   Relative distance to minimum upper bound &       &       &       &       &       &       &       &       &       &       &       \\
           mean 	& -0.06	&0.25	&0.13	&0.04	&0.03	&0.21	&0.07	&0.06	&0.04	&0.31 &0.11
\\
           std err & 0.01	&0.04	&0.02	&0.01	&0.01	&0.03	&0.02	&0.01	&0.01	&0.04 &0.01
\\
           0.25 q	&  -0.10	&0.05	&0.04	&0.00	&0.00	&0.06	&0.00	&0.01	&0.00	&0.10 &0.03
\\
           0.5 q 	&  -0.04	&0.10	&0.09	&0.01	&0.01	&0.12	&0.03	&0.03	&0.01	&0.16 &0.06
\\
           0.75 q 	&  -0.01	&0.23	&0.17	&0.05	&0.04	&0.31	&0.08	&0.06	&0.05	&0.34 &0.17
\\
\SetCell[c=12]{l}   Relative distance to maximum lower bound &       &       &       &       &       &       &       &       &       &       &       \\
           mean  	& -0.11	&0.64	&0.19	&0.01	&0.01	&0.23	&0.18	&0.03	&0.02	&0.73 &0.18
\\
           std err & 0.01	&0.03	&0.01	&0.00	&0.00	&0.01	&0.06	&0.00	&0.00	&0.02 &0.01
\\
           0.25 q 	& -0.17	&0.58	&0.10	&0.00	&0.00	&0.15	&0.00	&0.01	&0.00	&0.65 &0.09
\\
           0.5 q 	& -0.08	&0.69	&0.19	&0.01	&0.00	&0.21	&0.02	&0.02	&0.01	&0.74 &0.19
\\
           0.75 q 	& -0.05	&0.80	&0.29	&0.02	&0.01	&0.30	&0.37	&0.05	&0.02	&0.88 &0.26
\\
\hline
\end{tblr}
\end{adjustbox}
\caption{Relative distances yielded by computations in Block B with selected methods using CPLEX. The minima and the maximum are calculated over the Benders methods.}
\label{tab:rel distances block B CPLEX}
\end{table}

\section{Conclusion}
\label{sec:conclusion}

Our aim was to produce experimental evidence about the computational properties of a selection of methods relying on the latest variants of Benders cuts for the solution of deterministic MILPs. We presented an extensive list of unified and non-unified Benders cuts, expressing under a common mathematical structure the CGPs, Benders cuts, and Benders cycle, with a level of detail sufficient to conduct applications. We explained why diversified sets of MCFNDP instances constitute proper testing grounds and detailed how the formulation of the generic inceptive problem is specialized to conform to the specifications of the MCFNDP. We supplied element-explicit descriptions for every CGP and Benders cut considered.

We performed a systematic empirical analysis comparing performances from a broad selection of computational methods featuring unified and/or non-unified Benders cuts that also included direct branch-and-cut solutions with a MILP solver. The testing bench was made up of the complete set of Canad R instances of MCFNDP and assessments were based on bespoke evaluation methods as well as performance profiles. Except for the methods featuring the Glomb et al. cut, static inclusion of the strong capacity constraints was everywhere beneficial. The analysis carried with Gurobi as the supporting solver of the SMS++ computational library identified four best-performing Benders methods: methods involving the static version of the Brandenberg-Stursberg cut without and with inversion of objective and normalization constraint (F\_U\_StatBS and F\_StatBS) and methods involving the Hosseini-Turner cut equipped with the $l_1$ norm or the Hosseini-Turner bound on the $l_1$ norm (F\_R-1\_HT, F\_R-2\_HT). In the objective of comparing the computational characteristics of Gurobi and CPLEX and to verify their agreement on a set of leading methods, we also performed experiments with CPLEX as the supporting solver. We found out that CPLEX and Gurobi lead to similar orderings between methods, that CPLEX reached erroneous solutions for several testing instances when performing branch-and-Benders-cut, and that CPLEX's speed was lesser when performing branch-and-Benders-cut.

We envision the following extensions: \emph{First}, expanding the current testing ground by including instances of MCFNDP produced with the generator of \cite{LarsEtAl2023} so as to challenge and possibly reinforce our current conclusions. The generated instances would stem from alternative topologies and, once again, would cover a broad range in size and complexity. \emph{Second}, performing a parallel set of experiments where a first-tier open-source solver would support SMS++. The SCIP \citep{SCIP2026} MILP solver adjoined with HiGHS \citep{HiGH2026} as the auxiliary LP solver would be a natural candidate. \emph{Third}, conducting experiments with a separate testing bench consisting of instances originating from a differing MILP architecture. This would make it possible to probe the reliability of extraneous inferences by assessing the agreement of conclusions reached independently from the separate sets of experiments. \emph{Fourth}, it would be worthwhile to extend the comparisons performed in this paper to a stochastic setting. This extension is important in view of the large and often prohibitive resources required by stochastic experiments of realistic scales.

\section*{Acknowledgments}

This research was funded by the Canadian National Railway Company (CN) Chair in Optimization of Railway Operations at Université de Montréal, the Canada Research Chair program [950-232244] and IVADO fundamental research program on integrated machine learning and optimization for decision making under uncertainty. We are grateful to Enrico Calandrini and Donato Meoli for their expert help in improving our usage of the SMS++ library and in implementing new functionalities. We are also grateful to Serge Bisaillon for building the first prototype of the program.

\section*{Appendix: Element-explicit formulations} 
\setcounter{subsection}{0}
\renewcommand{\thesubsection}{A.\arabic{subsection}}

We express the CGPs and the corresponding cuts of Sections~\ref{sec:Benders algo} and \ref{sec:unified cuts}, as specialized according to Section~\ref{sec:matrix representation MCFNDP} to instantiate the MCFNDP in terms of the individual elements of their constituent vectors and matrices. These formulations are required when aiming at higher-speed calculations performed with a lower level programming language such as C++. Presence of the variables $\hat{\pi}^{(3)}$ and $\pi^{(3)}$ reflects the inclusion of all strong capacity constraints in the mathematical representation of the MCFNDP. Exclusion of some strong capacity constraints can be effected (i) by removing from objective and constraints every term where a corresponding coordinate of the $\hat{\pi}^{(3)}$ and $\pi^{(3)}$ appears.

\subsection{Standard feasibility and optimality cuts}

\paragraph{DSP:}
\begin{align}
  \max_{\hat{\pi}^{(1)}, \hat{\pi}^{(2)}, \hat{\pi}^{(3)}}\; & -\sum_i \sum_k w_i^k \hat{\pi}_{ik}^{(1)}
  - \sum_{ij}  u_{ij} \bar{y}_{ij} \hat{\pi}_{ij}^{(2)}
  - \sum_{ij} \sum_{k}  b_{ij}^k \bar{y}_{ij} \hat{\pi}_{ij}^{(3)k} \\
  \text{subject to }\; & I(j \in \mathcal{N}_i^+) \hat{\pi}_{ik}^{(1)} - I(i \in \mathcal{N}_j^-) \hat{\pi}_{jk}^{(1)}
  + \hat{\pi}_{ij}^{(2)} + \hat{\pi}_{ij}^{(3)k} + c_{ij}^k \geq 0\; \forall ij, \forall k, \\
   & \hat{\pi}^{(2)}, \hat{\pi}^{(3)} \geq 0.  \\
\end{align}

\paragraph{Feasibility cut (if DSP is unbounded):}
\begin{equation}  
  \sum_{ij} \hat{\pi}_{ij}^{(2)} u_{ij} y_{ij}
   + \sum_{ij} \sum_{k} \hat{\pi}_{ij}^{(3)k} b_{ij}^k y_{ij}
  \geq - \sum_i \sum_k \hat{\pi}_{ik}^{(1)} w_i^k,
\end{equation}
where $(\hat{\pi}^{(1)}, \hat{\pi}^{(2)}, \hat{\pi}^{(3)})$ is a ray of DSP.

\paragraph{Optimality cut (if DSP is bounded):}
\begin{equation}  
  \sum_{ij} \hat{\pi}_{ij}^{(2)} \, u_{ij} y_{ij}
   + \sum_{ij} \sum_{k} \hat{\pi}_{ij}^{(3)k} \, b_{ij}^k y_{ij}
  + \eta \geq - \sum_i \sum_k \hat{\pi}_{ik}^{(1)} \, w_i^k,
\end{equation}
where $(\hat{\pi}^{(1)}, \hat{\pi}^{(2)}, \hat{\pi}^{(3)})$ is an optimal solution of DSP.

\subsection{Optimality cut of Magnanti-Wong}

\paragraph{CGP-M:}
\begin{align}
  \max_{\hat{\pi}^{(1)}, \hat{\pi}^{(2)}, \hat{\pi}^{(3)}}\; & -\sum_i \sum_k w_i^k \hat{\pi}_{ik}^{(1)}
  - \sum_{ij}  u_{ij} y^*_{ij} \hat{\pi}_{ij}^{(2)}
  - \sum_{ij} \sum_{k}  b_{ij}^k y^*_{ij} \hat{\pi}_{ij}^{(3)k} \\
  \text{subject to }\; & I(j \in \mathcal{N}_i^+) \hat{\pi}_{ik}^{(1)} - I(i \in \mathcal{N}_j^-) \hat{\pi}_{jk}^{(1)}
  + \hat{\pi}_{ij}^{(2)} + \hat{\pi}_{ij}^{(3)k} + c_{ij}^k \geq 0\; \forall ij, \forall k, \\
& -\sum_i \sum_k w_i^k \hat{\pi}_{ik}^{(1)}
  - \sum_{ij}  u_{ij} \bar{y}_{ij} \hat{\pi}_{ij}^{(2)}
  - \sum_{ij} \sum_{k}  b_{ij}^k \bar{y}_{ij} \hat{\pi}_{ij}^{(3)k} = \widehat{\textit{V}}(\bar{y}), \\  
   & \hat{\pi}^{(2)}, \hat{\pi}^{(3)} \geq 0.
\end{align}

\paragraph{Optimality cut:}
\begin{equation}  
  \sum_{ij} \hat{\pi}_{ij}^{(2)} \, u_{ij} y_{ij}
   + \sum_{ij} \sum_{k} \hat{\pi}_{ij}^{(3)k} \, b_{ij}^k y_{ij}
  + \eta \geq - \sum_i \sum_k \hat{\pi}_{ik}^{(1)} \, w_i^k,
\end{equation}
where $(\hat{\pi}^{(1)}, \hat{\pi}^{(2)}, \hat{\pi}^{(3)})$ is an optimal solution of CGP-M.

\subsection{Optimality cut of Papadakos}

\paragraph{CGP-P:}
\begin{align}
  \max_{\hat{\pi}^{(1)}, \hat{\pi}^{(2)}, \hat{\pi}^{(3)}}\; & -\sum_i \sum_k w_i^k \hat{\pi}_{ik}^{(1)}
  - \sum_{ij}  u_{ij} y^*_{ij} \hat{\pi}_{ij}^{(2)}
  - \sum_{ij} \sum_{k}  b_{ij}^k y^*_{ij} \hat{\pi}_{ij}^{(3)k} \\
  \text{subject to }\; & I(j \in \mathcal{N}_i^+) \hat{\pi}_{ik}^{(1)} - I(i \in \mathcal{N}_j^-) \hat{\pi}_{jk}^{(1)}
  + \hat{\pi}_{ij}^{(2)} + \hat{\pi}_{ij}^{(3)k} + c_{ij}^k \geq 0\; \forall ij, \forall k, \\
   & \hat{\pi}^{(2)}, \hat{\pi}^{(3)} \geq 0.
\end{align}

\paragraph{Optimality cut:}
\begin{equation}  
  \sum_{ij} \hat{\pi}_{ij}^{(2)} \, u_{ij} y_{ij}
   + \sum_{ij} \sum_{k} \hat{\pi}_{ij}^{(3)k} \, b_{ij}^k y_{ij}
  + \eta \geq - \sum_i \sum_k \hat{\pi}_{ik}^{(1)} \, w_i^k,
\end{equation}
where $(\hat{\pi}^{(1)}, \hat{\pi}^{(2)}, \hat{\pi}^{(3)})$ is an optimal solution of CGP-P.

\subsection{Unified cuts of Fischetti et al., Brandenberg-Stursberg and Conforti-Wolsey}
\label{sec_appendix: unified cuts of FBC}

\paragraph{CGP-FBC:}
\begin{align}
  \max_{\hat{\pi}^{(1)}, \hat{\pi}^{(2)}, \hat{\pi}^{(3)}, \hat{\pi}_0}\; & \Phi \\
  \text{subject to }\; & \sum_i \sum_k w_i^k \hat{\pi}_{ik}^{(1)}
  + \sum_{ij} u_{ij} \bar{y}_{ij} \hat{\pi}_{ij}^{(2)}
  + \sum_{ij} \sum_{k} b_{ij}^k \bar{y}_{ij} \hat{\pi}_{ij}^{(3)k}
  + \bar{\eta} \hat{\pi}_0 \leq -1 \\
  &  I(j \in \mathcal{N}_i^+) \hat{\pi}_{ik}^{(1)} - I(i \in \mathcal{N}_j^-) \hat{\pi}_{jk}^{(1)}
  + \hat{\pi}_{ij}^{(2)} + \hat{\pi}_{ij}^{(3)k} + c_{ij}^k \hat{\pi}_0 \geq 0\; 
  \forall ij,\, \forall k, \\
  & \hat{\pi}^{(2)}, \hat{\pi}^{(3)}, \hat{\pi}_0 \geq 0, 
\end{align}
where
\[
\Phi :\equiv
\begin{cases}
-\sum_{ij} \sum_k \hat{\pi}_{ij}^{(3)k} - \sum_{ij} \hat{\pi}_{ij}^{(2)} - \hat{\pi}_0, & \text{if Fischetti et al. cut}, \\
-\sum_{ij} \sum_k b_{ij}^k \hat{\pi}_{ij}^{(3)k} - \sum_{ij} u_{ij} \hat{\pi}_{ij}^{(2)} - \hat{\pi}_0,
& \text{if static Brandenberg-Stursberg cut}, \\
-\sum_k  b_{ij}^k (\tilde{y}_{ij} - \bar{y}_{ij}) \hat{\pi}_{ij}^{(3)k},\\
  - \sum_{ij} u_{ij} (\tilde{y}_{ij} - \bar{y}_{ij}) \hat{\pi}_{ij}^{(2)}
  - (\tilde{\eta} - \bar{\eta}) \hat{\pi}_0, & \text{if adaptive Brandenberg-Stursberg cut}.
\end{cases}
\]

\paragraph{CGP-FBC-Invert:}
\begin{align}
  \max_{\hat{\pi}^{(1)}, \hat{\pi}^{(2)}, \hat{\pi}^{(3)}, \hat{\pi}_0}\; & -\sum_i \sum_k w_i^k \hat{\pi}_{ik}^{(1)}
  - \sum_{ij} u_{ij} \bar{y}_{ij} \hat{\pi}_{ij}^{(2)}
  - \sum_{ij} \sum_{k} b_{ij}^k \bar{y}_{ij} \hat{\pi}_{ij}^{(3)k}
  - \bar{\eta} \hat{\pi}_0 \\
  \text{subject to }\; & I(j \in \mathcal{N}_i^+) \hat{\pi}_{ik}^{(1)} - I(i \in \mathcal{N}_j^-) \hat{\pi}_{jk}^{(1)}
  + \hat{\pi}_{ij}^{(2)} + \hat{\pi}_{ij}^{(3)k} + c_{ij}^k \hat{\pi}_0 \geq 0\; \forall ij,\; \forall k, \\
  & \Psi, \\
   & \hat{\pi}^{(2)}, \hat{\pi}^{(3)}, \hat{\pi}_0 \geq 0,
\end{align}
where
\[
\Psi :\equiv
\begin{cases}
-\sum_{ij} \sum_k \hat{\pi}_{ij}^{(3)k} - \sum_{ij} \hat{\pi}_{ij}^{(2)} - \hat{\pi}_0 = -1, & \text{if Fischetti et al. cut}, \\
-\sum_{ij} \sum_k b_{ij}^k \hat{\pi}_{ij}^{(3)k} - \sum_{ij} u_{ij} \hat{\pi}_{ij}^{(2)} - \hat{\pi}_0 = -1, & \text{if static Brandenberg-Stursberg cut}, \\
-\sum_{ij} \sum_k b_{ij}^k (\tilde{y}_{ij} - \bar{y}_{ij}) \hat{\pi}_{ij}^{(3)k} \\
  - \sum_{ij} u_{ij} (\tilde{y}_{ij} - \bar{y}_{ij}) \hat{\pi}_{ij}^{(2)}
  - (\tilde{\eta} - \bar{\eta})\, \hat{\pi}_0 = -1, & \text{if adaptive Brandenberg-Stursberg cut}.
\end{cases}
\]

\paragraph{Unified cut:}
\begin{equation}  
  \sum_{ij} \hat{\pi}_{ij}^{(2)} \, u_{ij} y_{ij}
   + \sum_{ij} \sum_{k} \hat{\pi}_{ij}^{(3)k} \, b_{ij}^k y_{ij}
  + \hat{\pi}_0 \eta \geq - \sum_i \sum_k \hat{\pi}_{ik}^{(1)} \, w_i^k,
 \end{equation}
where $(\hat{\pi}^{(1)}, \hat{\pi}^{(2)}, \hat{\pi}^{(3)}, \hat{\pi}_0)$ is an optimal solution of CGP-FBC or CGP-FBC-Invert.

\subsection{Unified cut of Seo et al.}

\paragraph{CGP-S:}
\begin{align}
  \max_{\hat{\pi}^{(1)}, \hat{\pi}^{(2)}, \hat{\pi}^{(3)}, \hat{\pi}_0}\; & -\sum_i \sum_k w_i^k \hat{\pi}_{ik}^{(1)}
  - \sum_{ij}  u_{ij} \tilde{y}_{ij} \hat{\pi}_{ij}^{(2)}
  - \sum_{ij} \sum_{k}  b_{ij}^k \tilde{y}_{ij} \hat{\pi}_{ij}^{(3)k}
  - \tilde{\eta} \hat{\pi}_0 \\
  \text{subject to }\; & \sum_{ij} \sum_k  b_{ij}^k (\tilde{y}_{ij} - \bar{y}_{ij}) \hat{\pi}_{ij}^{(3)k}
  + \sum_{ij}  u_{ij} (\tilde{y}_{ij} - \bar{y}_{ij}) \hat{\pi}_{ij}^{(2)}
  + (\tilde{\eta} - \bar{\eta}) \hat{\pi}_0 = 1, \\
   & I(j \in \mathcal{N}_i^+) \hat{\pi}_{ik}^{(1)} - I(i \in \mathcal{N}_j^-) \hat{\pi}_{jk}^{(1)}
  + \hat{\pi}_{ij}^{(2)} + \hat{\pi}_{ij}^{(3)k} + c_{ij}^k \hat{\pi}_0 \geq 0\; \forall ij, \forall k, \\
   & \hat{\pi}^{(2)}, \hat{\pi}^{(3)}, \hat{\pi}_0 \geq 0.
\end{align}

\paragraph{Unified cut: (formally identical to that of Section~\ref{sec_appendix: unified cuts of FBC})}
\begin{equation}  
  \sum_{ij} \hat{\pi}_{ij}^{(2)} \, u_{ij} y_{ij}
   + \sum_{ij} \sum_{k} \hat{\pi}_{ij}^{(3)k} \, b_{ij}^k y_{ij}
  + \hat{\pi}_0 \eta \geq - \sum_i \sum_k \hat{\pi}_{ik}^{(1)} \, w_i^k,
\end{equation}
where $(\hat{\pi}^{(1)}, \hat{\pi}^{(2)}, \hat{\pi}^{(3)}, \hat{\pi}_0)$ is an optimal solution of CGP-S.

\subsection{Unified cut of Hosseini-Turner}

\paragraph{CGP-H:}
\begin{align}
  \max_{\pi^{(1)}, \pi^{(2)}, \pi^{(3)}, \pi_0}\; & -\sum_i \sum_k w_i^k \pi_{ik}^{(1)}
  - \sum_{ij}  u_{ij} \bar{y}_{ij} \pi_{ij}^{(2)}
  - \sum_{ij} \sum_{k}  b_{ij}^k \bar{y}_{ij} \pi_{ij}^{(3)k}
  - (\bar{\theta} - \sum_{ij} f_{ij} \bar{y}_{ij}) \pi_0 \\
  \text{subject to }\; & I(j \in \mathcal{N}_i^+) \pi_{ik}^{(1)} - I(i \in \mathcal{N}_j^-) \pi_{jk}^{(1)}
  + \pi_{ij}^{(2)} + \pi_{ij}^{(3)k} + c_{ij}^k \pi_0 \geq 0\; \forall ij,\, \forall k, \\
  & \Gamma, \\
   & \pi^{(2)}, \pi^{(3)}, \pi_0 \geq 0,
\end{align}
where the set of constraints $\Gamma$ is defined as follows according to the choice of an $l_p$ norm or its relaxation: 
\paragraph{$l_p, p = \infty, q = 1$}
\[
\Gamma :\equiv
\begin{cases}
f_{ij} \pi_0 - \sum_k b_{ij}^k \pi_{ij}^{(3)k} - u_{ij} \pi_{ij}^{(2)} \geq -1\; \forall{ij}, \\
f_{ij} \pi_0 - \sum_k b_{ij}^k \pi_{ij}^{(3)k} - u_{ij} \pi_{ij}^{(2)} \leq 1\; \forall{ij}, \\
\pi_0 \leq 1.
\end{cases}
\]

\paragraph{$l_p, p = 1, q = \infty$}
\[
\Gamma :\equiv
\begin{cases}
\pi_0 + \sum_{ij} \tau_{ij} \leq 1, \\
f_{ij} \pi_0 - \sum_k b_{ij}^k \pi_{ij}^{(3)k} - u_{ij} \pi_{ij}^{(2)} \geq -\tau_{ij}\; \forall{ij}, \\
f_{ij} \pi_0 - \sum_k b_{ij}^k \pi_{ij}^{(3)k} - u_{ij} \pi_{ij}^{(2)} \leq \tau_{ij}\; \forall{ij}, \\
\tau_{ij} \geq 0\; \forall{ij}.
\end{cases}
\]

\paragraph{relaxation of $l_p, p=1$}
\[
\Gamma :\equiv (1 + \sum_{ij} f_{ij}) \pi_0 + \sum_{ij} \sum_k b_{ij}^k \pi_{ij}^{(3)k} + \sum_{ij} u_{ij} \pi_{ij}^{(2)} \leq 1.
\]

\paragraph{Unified cut:}
\begin{equation}  
  \sum_{ij} \pi_{ij}^{(2)} \, u_{ij} y_{ij}  
   + \sum_{ij} \sum_{k} \pi_{ij}^{(3)k} \, b_{ij}^k y_{ij}
   - \sum_{ij} \pi_0 f_{ij} y_{ij}   
  + \pi_0 \theta  
  \geq - \sum_i \sum_k \pi_{ik}^{(1)} \, w_i^k,
\end{equation}
where $(\pi^{(1)}, \pi^{(2)}, \pi^{(3)}, \pi_0)$ is an optimal solution of CGP-H.

\subsection{Unified cut of Glomb et al.}

\paragraph{CGP-G:}
\begin{align}
  \max_{\pi^{(1)}, \pi^{(2)}, \pi^{(3)}, \pi_0}\; &
  -\sum_{ij} \sum_k b_{ij}^k \bar{y}_{ij} \pi_{ij}^{(3)k} - \sum_{ij}  u_{ij} \bar{y}_{ij} \pi_{ij}^{(2)} - \sum_i \sum_k w_i^k \pi_{ik}^{(1)} + \Big(\sum_{ij} f_{ij} \bar{y}_{ij} - \beta\Big) \pi_0 \\
  \text{subject to }\; & -\sum_{ij} \sum_k  b_{ij}^k\, (\bar{y}_{ij} - \tilde{y}_{ij}) \pi_{ij}^{(3)k} - \sum_{ij} u_{ij} (\bar{y}_{ij} - \tilde{y}_{ij}) \pi_{ij}^{(2)}
  + \sum_{ij} f_{ij}(\bar{y}_{ij} - \tilde{y}_{ij}) \pi_0 = 1, \\
& I(j \in \mathcal{N}_i^+) \pi_{ik}^{(1)} - I(i \in \mathcal{N}_j^-) \pi_{jk}^{(1)}
  + \pi_{ij}^{(2)} + \pi_{ij}^{(3)k} + c_{ij}^k \pi_0 \geq 0\; \forall ij,\, \forall k, \\
  & \pi^{(3)}, \pi^{(2)}, \pi_0 \geq 0 
\end{align}

\paragraph{Unified cut:}
\begin{equation}
  -\sum_{ij} \sum_k \pi_{ij}^{(3)k}\, b_{ij}^k y_{ij}
  - \sum_{ij} \pi_{ij}^{(2)}\, u_{ij} y_{ij} - \pi_0 \eta
  \leq \sum_i \sum_k \pi_{ik}^{(1)} \, w_i^k
  \label{eq:10.25}
\end{equation}
where $(\pi^{(1)}, \pi^{(2)}, \pi^{(3)}, \pi_0)$ is an optimal solution of CGP-G.

\paragraph{DSP-G:}
\begin{align}
  \max_{\hat{\pi}^{(1)}, \hat{\pi}^{(2)}, \hat{\pi}^{(3)}}\; & -\sum_i \sum_k w_i^k \hat{\pi}_{ik}^{(1)}
  - \sum_{ij}  u_{ij} \bar{y}_{ij} \hat{\pi}_{ij}^{(2)}
  - \sum_{ij} \sum_{k}  b_{ij}^k \bar{y}_{ij} \hat{\pi}_{ij}^{(3)k} \\
  \text{subject to }\; & I(j \in \mathcal{N}_i^+) \hat{\pi}_{ik}^{(1)} - I(i \in \mathcal{N}_j^-) \hat{\pi}_{jk}^{(1)}
  + \hat{\pi}_{ij}^{(2)} + \hat{\pi}_{ij}^{(3)k} + c_{ij}^k \geq 0\; \forall ij, \forall k, \\
&-\sum_i \sum_k w_i^k \hat{\pi}_{ik}^{(1)}
  - \sum_{ij}  u_{ij} \tilde{y}_{ij} \hat{\pi}_{ij}^{(2)}
  - \sum_{ij} \sum_{k}  b_{ij}^k \tilde{y}_{ij} \hat{\pi}_{ij}^{(3)k} = \beta - \sum_{ij} f_{ij} \tilde{y}_{ij}, \\  
   & \hat{\pi}^{(2)}, \hat{\pi}^{(3)} \geq 0.
\end{align}

\bibliographystyle{plainnat_custom}

\bibliography{refs}

\end{document}